\documentclass[a4paper,11pt]{article}
\usepackage{amsmath, amssymb, amsthm}
\usepackage{graphicx} 

\usepackage[latin1]{inputenc}

\begin{document}

\def\Fiv{{\rm Div}}\def\N{{\mathbb N}} \def\Div{{\rm Div}}
\def\Z{{\mathbb Z}} \def\Ker{{\rm Ker}} \def\PP{{\mathcal P}}
\def\CC{{\mathcal C}} \def\LL{{\mathcal L}} \def\RR{{\mathcal R}}
\def\JJ{{\mathcal J}} \def\Card{{\rm Card}} \def\Inv{{\rm Inv}}
\def\Id{{\rm Id}} \def\T{{\mathbb T}} \def\si{{\rm si}}
\def\hgt{{\rm hgt}} \def\AA{{\mathcal A}}\def\SS{{\mathcal S}}


\title{\bf{PreGarside monoids and groups, parabolicity, amalgamation, and FC property}}

\author{\textsc{Eddy Godelle\footnotemark[1] and Luis Paris\footnote{Both authors are partially supported by the {\it Agence Nationale de la Recherche} ({\it projet Théorie de Garside}, ANR-08-BLAN-0269-03).}}}

\date{\today}

\maketitle

\begin{abstract}
\noindent
We define the notion of preGarside group slightly lightening the definition of Garside group so that all Artin-Tits groups are preGarside groups. This paper intends to give a first basic study on these groups. Firstly, we introduce the notion of parabolic subgroup, we prove that any preGarside group has a (partial) complemented presentation, and we characterize the parbolic subgroups in terms of these presentations. Afterwards we prove that the amalgamated product of two preGarside groups along a common parabolic subgroup is again a preGarside group. This enables us to define the family of preGarside groups of FC type as the smallest family of preGarside groups that contains the Garside groups and that is closed by amalgamation along parabolic subgroups. Finally, we make an algebraic and combinatorial study on FC type preGarside groups and their parabolic subgroups. 
\end{abstract}


\noindent
{\bf AMS Subject Classification.} Primary: 20F36.  


\section{Introduction}

Let $S$ be a (non necessarily finite) set. A \emph{Coxeter matrix} over $S$ is a square matrix $M = (m_{s,t})_{s,t \in S}$ indexed by the elements of $S$ and satisfying (1) $m_{s,s}=1$ for all $s \in S$; (2) $m_{s,t} = m_{t,s} \in \{2,3,4, \dots, \infty\}$ for all $s,t \in S$, $s \neq t$. With this combinatorial data one can associate an \emph{Artin-Tits group}, which we denote by $G=G_M$, and which is combinatorially defined by the generating set $S$ and the relations
\[
\underbrace{sts \cdots}_{m_{s,t}\ \mathrm{terms}} = \underbrace{tst \cdots}_{m_{s,t}\ \mathrm{terms}}\,, \quad \text{for } s,t \in S\,,\ s \neq t \text{ and } m_{s,t} \neq \infty\,.
\]
The \emph{Coxeter group} associated with $M$, denoted by $W=W_M$, is the quotient of $G$ by the relations $s^2=1$, $s \in S$.

\bigskip\noindent
An Artin-Tits group $G_M$ is called of {\it spherical type} if the Coxeter group $W_M$ of $M$ is finite. In the early 1970s, inspired in particular by the work of Arnold \cite{Arnol1, Arnol2, Arnol3, Arnol4} on the cohomolgy of braid groups and the work of Garside \cite{Garsi1} on the conjugacy problem of braid groups, Brieskorn  \cite{Bries1, Bries2, Bries3}, Saito \cite{BriSai1}, and Deligne \cite{Delig1} (see also \cite{Lin1, Gorju1}) initiated the study of spherical type Artin-Tits groups  as well as they links with hyperplane arrangements. These groups are now very well-understood. In particular, they are known to have solvable word and conjugacy problems, and to be biautomatic \cite{Charn1, Charn2}, and the spaces of regular orbits of the associated Coxeter groups are classifying spaces for them \cite{Delig1}.

\bigskip\noindent
The next important step in the study of Artin-Tits groups was Van der Lek's thesis \cite{VdL} whose main result is that every Artin-Tits group is the fundamental group of the space of regular orbits of the associated Coxeter group acting on the complexified Tits cone. But, it also contains a study on parabolic subgroups of Artin-Tits groups, as well as the statement of the so-called $K(\pi,1)$ conjecture for Artin-Tits groups, one of the central questions in the subject. Recall that a \emph{standard parabolic subgroup} of $G_M$ is defined to be a subgroup generated by a subset of $S$, and a \emph{parabolic subgroup} is a subgroup conjugate to a standard parabolic subgroup. By \cite{VdL} (see also \cite{Paris1, GodPar1}), a standard parabolic subgroup is itself an Artin-Tits group in a canonical way. 

\bigskip\noindent
In \cite{ChaDav1} Charney and Davis used techniques from hyperbolic geometry and geometric group theory (CAT(0) spaces) to solve the $K(\pi,1)$ conjecture on two new families of Artin-Tits groups: that of Artin-Tits groups of FC type, and that of Artin-Tits groups of dimension 2 (see also \cite{ChaDav2}). The family of Artin-Tits groups of \emph{FC type} is the smallest family of Artin-Tits groups that contains the Artin-Tits groups of spherical type and that is closed under amalgamation over standard parabolic subgroups. On the other hand, an Artin-Tits group $G_M$ is of \emph{dimension $2$} if, for every subset $X$ of $S$ of cardinality at least $3$, the parabolic subgroup of $G_M$ generated by $X$ is not of spherical type. The word problem is known to be solvable for these groups \cite{Cherma1, Alt, AltCha1}, but it is not known whether they have solvable conjugacy problem. An algebraic and combinatorial study of parabolic subgroups of these groups can be found in \cite{Godel1, Godel2}.

\bigskip\noindent
Two notions play a prominent role in all these studies: that of parabolic subgroup (already defined), and that  of Artin-Tits monoid. The \emph{Artin-Tits monoid} associated with the Coxeter matrix $M$ is the monoid having as monoid presentation the same presentation as $G_M$ viewed as a group. By \cite{Par}, this embeds into $G_M$.

\bigskip\noindent
Inspired by Garside's work \cite{Garsi1} and Thurston's work \cite{ECHLPT1} on braid groups, both extended to spherical type Artin-Tits groups (see \cite{Charn1, Charn2}), Dehornoy and the second author \cite{DePa} introduced in 1999 the notions of Garside monoids and Garside groups (see also \cite{Deh}), and they showed that these monoids and groups share many properties with Artin-Tits monoids and groups of spherical type such as solvable word and conjugacy problems, torsion freeness, and biautomaticity. Since then, Garside groups have become popular objects of study. Their definitions are given in Section 2. 

\bigskip\noindent
We define (see Section 2) the notion of preGarside monoid slightly lightening the definition of Garside monoid so that all Artin-Tits monoids are preGarside monoids. A preGarside group is defined to be the enveloping group of a preGarside monoid.  Although these notions are not new (see \cite{DiM, Deh3, Deh4}) (preGarside monoids are often called locally Garside monoids), there are no studies dedicated to these monoids and groups. Hence, the present paper may be considered as a first step to their study.

\bigskip\noindent
In Section 2 we define the notions of standard parabolic submonoid of a preGarside monoid and of standard parabolic subgroup of a preGarside group. These definitions extend, in the one hand, the notion of standard parabolic subgroup (resp. submonoid) of an Artin-Tits group (resp. monoid), and, in the other hand, the notion of standard parabolic subgroup (resp. submonoid) of a Garside group (resp. monoid)  \cite{God, Godel3}. We prove that any preGarside monoid (or group) has a (partial) complemented presentation in the style of the complemented presentations for Garside groups and monoids given in \cite{DePa} (see Theorem 2.6). Moreover, we show that a standard parabolic subgroup (resp. submonoid) is necessarily generated by some subset of the generating family of the (partial) complemented presentation, and give necessary and sufficent conditions for a subset of this generating family to span a standard parabolic subgroup (resp. submonoid) (see Theorem 2.8). 

\bigskip\noindent
The most significant result in Section 3 is that the amalgamated product $M_1\!*_N\!M_2$ of two preGarside monoids $M_1,M_2$ along a common standard parabolic submonoid $N$ is again a preGarside monoid (see Proposition 3.11). But, our study does not end with this result. Indeed, we also prove that the two monoids $M_1, M_2$ embed into $M_1\!*_N\!M_2$ (Proposition 3.1) --this is not true in general-- and that the amalgamated product $M_1\!*_N\!M_2$ admits normal forms similar to the standard normal forms for amalgamated products of groups (Proposition 3.3). Moreover, we characterize the standard parabolic submonoids of $M_1\!*_N\!M_2$ in terms of standard parabolic submonoids of $M_1$ and of $M_2$ (Proposition~3.12).

\bigskip\noindent
The fact that the amalgamated product of two preGarside monoids along a common standard parabolic submonoid is still a preGarside monoid enables us to construct new examples of preGarside groups (and monoids). It also enables us to define the family of preGarside monoids of FC type as the smallest family of preGarside monoids that contains the Garside monoids and that is closed by amalgamation along standard parabolic submonoids. This extends the definition of Artin-Tits monoids (and groups) of FC type given above.

\bigskip\noindent
Section 4 is dedicated to the algebraic and combinatorial study of preGarside groups of FC type and their standard parabolic subgroups. In particular, we prove the following.

\bigskip\noindent
{\bf Theorem 4.10.}
{\it Let $M$ be a preGarside monoid of FC type, and let $G(M)$ be its enveloping group.
\begin{itemize}
\item[(P1)]
The natural morphism $\iota : M \to G(M)$ is injective.
\item[(P2)]
Let $N$ be a standard parabolic submonoid of $M$. The standard parabolic subgroup of $G(M)$ generated by $N$ is isomorphic to $G(N)$, and we have $G(N) \cap M = N$.
\item[(P3)]
Let $N, N'$ be standard parabolic submonoids of $M$. Then $N \cap N'$ is a standard parabolic submonoid, and $G(N) \cap G(N') = G(N \cap N')$.
\item[(P4)]
$G(M)$ is torsion free.
\end{itemize}}

\bigskip\noindent
Properties (P1), (P2), and (P3) of the above theorem are known to hold for all Artin-Tits monoids \cite{Par,VdL}.  Actually, the second part of Property (P2) is not proved (nor stated) in \cite{VdL}, but it can be easily deduced from it. However, to know whether Property (P4) holds for all Artin-Tits monoids is an open question.

\bigskip\noindent
Concerning the algorithmic properties of FC type preGarside groups we prove the following.

\bigskip\noindent
{\bf Corollary 4.16.}
{\it Let $M$ be a finitely generated preGarside monoid of FC type, let $S$ be its set of atoms (which generates $M$), and let $G(M)$ be the enveloping group of $M$.
\begin{itemize}
\item[(1)]
$G(M)$ has a solution to the word problem.
\item[(2)]
There exists an algorithm which, given $w \in S^{\pm *}$, decides whether the element $\overline{w}$ in $G(M)$ represented by $w$ belongs to $M$ or not.
\item[(3)]
Let $H$ be a standard parabolic subgroup of $G(M)$. There exists an algorithm which, given $w \in S^{\pm *}$, decides whether $\overline{w} \in H$.
\end{itemize}}


\section{Parabolic submonoids and subgroups, and presentations}

\subsection{Definitions and basic properties}

We start with some terminology. Consider a monoid $M$. It is said to be \emph{cancellative} if, for all $a,b,c,d \in M$, the equality $cad = cbd$ imposes $a = b$. An element $b$ is called a \emph{factor} of an element $a$ if we can write $a = cbd$ in $M$. We denote by $\Div(a)$ the set of factors of $a$. When $a = bc$, we say that $b$ \emph{left-divides} $a$ and write $b \preceq_L a$. Similarly, we say that $c$ \emph{right-divides} $a$ and write $c \preceq_R a$.  An element $a$ is said to be \emph{balanced} if its sets of right-divisors and of left-divisors are equal, which in this case have to be equal to $\Div(a)$. We say that $M$ is \emph{atomic} if there exists a mapping $\nu : M \to \N$, called a \emph{norm}, satisfying $\nu(a) >0$ for $a \neq 1$ and $\nu(ab) \ge \nu(a) + \nu(b)$ for all $a,b \in M$. Note that the existence of such a mapping implies that the relations $\preceq_L$ and $\preceq_R$ are partiel orders on $M$. 
The enveloping group of a monoid $M$ will be always denoted by $G(M)$, and the canonical morphism $M \to G(M)$ by $\iota = \iota_M : M \to G(M)$.

\bigskip\noindent
{\bf Definition.}
A monoid $M$ is said to be a \emph{preGarside monoid} if 
\begin{itemize}
\item[(a)] 
it is cancellative and atomic; 
\item[(b)] 
for all $a,b \in M$, if the set $\{c \in M\, |\, a \preceq_L c \text{ and } b \preceq_L c\}$ is nonempty, then it has a least element, denoted by $a \vee_{\!\scriptscriptstyle L} b$;
\item[(c)] 
for all $a,b \in M$, if the set $\{c \in M\, |\, a\preceq_R c \text{ and } b \preceq_R c\}$ is nonempty, then it has a least element, denoted by $a \vee_{\!\scriptscriptstyle R} b$.
\end{itemize}
A \emph{Garside element} of a preGarside monoid is a balanced element whose set of factors generates the whole monoid. When such an element exists, we say that the monoid is a \emph{Garside monoid}. A \emph{preGarside group} $G(M)$ is the enveloping group of a preGarside monoid $M$. Similarly, a \emph{Garside group} $G(M)$ is the enveloping group of a Garside monoid $M$.

\bigskip\noindent
{\bf Remark.}
We will not assume that our monoids are finitely generated, except when we will study algorithmic questions. Indeed, although our algorithmic results can be applied to some well-understood (pre)Garside monoids and groups such as the ones introduced by Digne in \cite{Digne1, Digne2}, a treatment of questions such as the word problem in the context of infinitely generated monoids and groups requires extra hypothesis such as a machine which recognizes the generating system of the given monoid or group. On the other hand, in most of the references (see \cite{DDGKM1, DePa} for instance), a requirement in the definition of a Garside monoid is that it is finitely generated, but this is not always true (see \cite{Digne1, Digne2}).
In this paper we remove this requirement, but the reader must understand that this is not completely standard.

\bigskip\noindent
As pointed out in the introduction, the seminal examples of Garside groups are the Artin-Tits groups of spherical type. Note also that all the Artin-Tits monoids are preGarside monoids, and hence all the Artin-Tits groups are preGarside groups (see \cite{BriSai1} and \cite{Mic}). We refer to \cite{Deh, DDGKM1} for the general theory on Garside groups.

\bigskip\noindent
Recall that an \emph{atom} in an atomic monoid $M$ is an element $a \in M$ satisfying $a=bc \Rightarrow$ $b=1$ or $c=1$ for all $b,c \in M$. We denote by $\AA(M)$ the set of atoms of $M$. Note that any generating set of $M$ contains $\AA(M)$. In particular, $M$ is finitely generated if and only if $\AA(M)$ is finite. 

\bigskip\noindent
{\bf Remark.} 
Let $\nu: M \to \N$ be a norm. Let $g \in M$. If $g = x_1 \cdots x_\ell$ is an expression of $g$ over the atoms, then $\ell \le \nu(g)$. In particular, the lengths of the expressions of $g$ over the atoms are bounded. Then, it is easily seen that the map $\tilde \nu : M \to \N$ which sends each $g \in M$ to the maximal length of an expression of $g$ over the atoms is a norm. 

\bigskip\noindent
{\bf Remark.}
A monoid $M$ is \emph{Noetherian} if every sequence $(a_n)_{n\in\N}$ of elements of $M$ such $a_{n+1}$ is a factor of $a_n$ stabilizes. It is easily seen that an atomic monoid is Noetherian, and, if $M$ is a finitely generated monoid, then $M$ is Noetherian if and only if it is atomic. Many of the results in the paper can be proved in the framework of Noetherian monoids, but the proofs are longer and more complicate. So, in order to simplify and shorten the proofs, we decide to make our study with atomic monoids.

\bigskip\noindent
{\bf Lemma 2.1.}
{\it Let $M$ be a preGarside monoid, and let $X \subset M$ be a nonempty subset.
\begin{itemize}
\item[(1)]
The set $\{a \in M\, |\, a \preceq_L x \text{ for all } x \in X\}$ has a greatest element (for the ordering $\preceq_L$), denoted by $\wedge_{\!\scriptscriptstyle L}X$. Similarly, the set $\{a \in M\, |\, a \preceq_R x \text{ for all } x \in X\}$ has a greatest element (for the ordering $\preceq_R$), denoted by $\wedge_{\!\scriptscriptstyle R}X$.
\item[(2)]
If the set $\{ a \in M\, |\, x \preceq_L a \text{ for all } x \in X\}$ is nonempty, then it has a least element  (for the ordering $\preceq_L$), denoted by $\vee_{\!\scriptscriptstyle L}X$. Similarly, if the set $\{ a \in M\, |\, x \preceq_R a \text{ for all } x \in X\}$ is nonempty, then it has a least element (for the ordering $\preceq_R$), denoted by $\vee_{\!\scriptscriptstyle R}X$.
\end{itemize}}

\bigskip\noindent
{\bf Proof.}
Let $Y = \{ y \in M \mid y \preceq_L x \text{ for all } x \in X\}$. Note that $1 \in Y$, thus $Y \neq \emptyset$. Let $\nu: M \to \N$ be a norm on $M$. The set $\{ \nu(y) \mid y \in Y\}$ is bounded by $\nu(x)$ for any $x \in X$, thus we may choose $y_0 \in Y$ such that $\nu(y_0)$ is maximal. If $y \in Y$, then $y \vee_{\!\scriptscriptstyle L} y_0$ exists and $( y \vee_{\!\scriptscriptstyle L} y_0) \preceq_L x$ for all $x \in X$, that is, $(y \vee_{\!\scriptscriptstyle L} y_0) \in Y$. Since $\nu(y_0)$ is maximal, it follows that $y_0 = (y \vee_{\!\scriptscriptstyle L} y_0)$, hence $y \preceq_L y_0$. So, $y_0 = \wedge_{\!\scriptscriptstyle L}X$. Now, set $Z=\{ z \in M \mid x \preceq_L z \text{ for all } x \in X\}$ and suppose $Z \neq \emptyset$. It is easily checked that $\wedge_{\!\scriptscriptstyle L}Z \in Z$, hence $\wedge_{\!\scriptscriptstyle L}Z$ is the least element of $Z$.
\qed

\bigskip\noindent
In the next proposition we gather some results on Garside monoids that we will need in the sequel. We refer to~\cite{Deh,DiM} for the proof.

\bigskip\noindent
{\bf Proposition 2.2.}
{\it Assume $M$ is a Garside monoid.
\begin{itemize}
\item[(1)]
The monoid $M$ has a (unique) minimal Garside element $\delta$, simply called {\rm the minimal Garside element} of $M$. 
\item[(2)]
$M$ is a lattice for left-divisibility and for right-divisibility. Furthermore, it injects into its enveloping group $G(M)$.
\item[(3)] 
Let $\Delta$ be a Garside element. Any element $a$ in $M$ has a unique decomposition $a_1\cdots a_n$ such that $a_n\neq 1$ and $a_i$ is the greatest element of $\Div(\Delta)$ that left-divides $a_i\cdots a_n$ for all $i \in \{1, \dots, n\}$.
\item[(4)]
Every element $g$ in $G(M)$ has a unique decomposition $ab^{-1}$ with $a,b$ in $M$ and $a\wedge_{\!\scriptscriptstyle R}b=1$.
\qed
\end{itemize}}

\bigskip\noindent
The decomposition in Proposition 2.2\,(3) is called \emph{left greedy normal form}. One can define a \emph{right greedy normal form} in a similar way. In this paper, by a \emph{greedy normal form} we will always mean a left greedy normal form. The decomposition in Proposition 2.2\,(4) is called \emph{right normal form}. One can also define a \emph{left normal form} in a similar way. From now on, by a \emph{normal form} we will always mean a right normal form. 

\subsection{Parabolic subgroups} 

In \cite{God} the first author introduced the notion of a standard parabolic subgroup of a Garside group. Here we extend this notion to the framework of preGarside groups.

\bigskip\noindent
{\bf Definition.}
Let $M$ be a monoid and let $N$ be a submonoid. We say that $N$ is \emph{special} if it is closed by factors, that is,  $ab\in N\implies a,b\in N$, for all $a,b \in M$.

\bigskip\noindent
{\bf Definition.}
Let $M$ be a preGarside monoid, and let $G(M)$ be its associated preGarside group. Denote by $\iota: M\to G(M)$ the canonical morphism. A submonoid $N$ of $M$ is said to be a \emph{standard parabolic submonoid} if 
\begin{itemize}
\item[(a)]
it is special;
\item[(b)]
for all $a,b \in N$, if $a\vee_{\!\scriptscriptstyle L} b$ exists, then $a \vee_{\!\scriptscriptstyle L} b \in N$, and if $a \vee_{\!\scriptscriptstyle R} b$ exists, then $a \vee_{\!\scriptscriptstyle R} b \in N$.
\end{itemize}
A standard parabolic submonoid is of \emph{spherical type} if it has a Garside element. A subgroup of $G(M)$ is a \emph{standard parabolic subgroup} if it is generated by the image $\iota(N)$ of a parabolic submonoid $N$ of $M$.
From now on, as we will never talk about general parabolic submonoids and subgroups, by a parabolic submonoid or subgroup we will mean a standard parabolic submonoid or subgroup.

\bigskip\noindent
{\bf Lemma 2.3.}
{\it Let $M$ be a preGarside monoid, and let $N$ be a parabolic submonoid of $M$.
\begin{itemize}
\item[(1)]
The monoid $N$ is a preGarside monoid. Moreover, it is a Garside monoid if and only if it is a spherical type submonoid of $M$.
\item[(2)]
Let $X$ be a non-empty subset of $N$. Then $\wedge_{\!\scriptscriptstyle L}X$ and $\wedge_{\!\scriptscriptstyle R}X$ belong to $N$. Similarly, $\vee_{\!\scriptscriptstyle L}X$ and $\vee_{\!\scriptscriptstyle R}X$ belong to $N$ when they exist.
\end{itemize}}

\bigskip\noindent
{\bf Proof.}
The only non-trivial part of the lemma is that $\vee_{\!\scriptscriptstyle L}X$ and $\vee_{\!\scriptscriptstyle R} X$ belong to $N$ when they exist. 
Suppose that $\vee_{\!\scriptscriptstyle L}X$ exists. 
Let $\nu: M \to \N$ be a norm.
Observe that, if $Y$ is a nonempty finite subset of $X$, then $\vee_{\!\scriptscriptstyle L}Y$ exists, $\vee_{\!\scriptscriptstyle L}Y \in N$, and $\vee_{\!\scriptscriptstyle L}Y \preceq_L \vee_{\!\scriptscriptstyle L}X$, thus $\nu(\vee_{\!\scriptscriptstyle L}Y) \le \nu(\vee_{\!\scriptscriptstyle L}X)$.
Now, choose a nonempty finite subset $Y_0$ of $X$ such that $\nu(\vee_{\!\scriptscriptstyle L}Y_0)$ is maximal. 
If there was $x \in X$ such that $x \not\preceq_L \vee_{\!\scriptscriptstyle L}Y_0$, then we would have $(\vee_{\!\scriptscriptstyle L} Y_0) \precneqq_L (\vee_{\!\scriptscriptstyle L}(Y_0 \cup\{x\}))$, thus $\nu (\vee_{\!\scriptscriptstyle L}Y_0) \lneq \nu(\vee_{\!\scriptscriptstyle L}(Y_0 \cup\{x\})$, which would contradict the maximality of $\nu(\vee_{\!\scriptscriptstyle L}Y_0)$. 
Hence, $x \preceq_L \vee_{\!\scriptscriptstyle L}Y_0$ for all $x \in X$, therefore $\vee_{\!\scriptscriptstyle L} X = \vee_{\!\scriptscriptstyle L} Y_0 \in N$.
\qed

\bigskip\noindent
As in the case of Artin-Tits groups, we can say more on parabolic submonoids when the preGarside monoid is a Garside monoid.

\bigskip\noindent
{\bf Lemma 2.4.}
{\it Let $M$ be a Garside monoid, let $\Delta$ be a Garside element of $M$, and let $N$ be a parabolic submonoid of $M$.
\begin{itemize}
\item[(1)]
The monoid $N$ is of spherical type. Moreover, there exists a Garside element $\Delta_N$ of $N$ such that $\Div(\Delta_N) = \Div(\Delta)\cap N$.  
\item[(2)]
An element of $G(N)$ has the same left-normal form (resp. right-normal form) in $G(N)$ as in $G(M)$.
\end{itemize}}

\bigskip\noindent
{\bf Proof.}
Every element of $\Div(\Delta) \cap N$ left-divides $\Delta$, thus $\vee_{\!\scriptscriptstyle L}(\Div(\Delta) \cap N)$ exists. 
It lies in $N$ by Lemma 2.3. 
Similarly, $\vee_{\!\scriptscriptstyle R}(\Div(\Delta) \cap N)$ exists and lies in $N$. 
Since $\vee_{\!\scriptscriptstyle L}(\Div(\Delta) \cap N)$ and $\vee_{\!\scriptscriptstyle R}(\Div(\Delta) \cap N)$ belong to $\Div(\Delta) \cap N$, we have $\vee_{\!\scriptscriptstyle L}(\Div(\Delta) \cap N) \preceq_R \vee_{\!\scriptscriptstyle R}(\Div(\Delta) \cap N)$ and $\vee_{\!\scriptscriptstyle R}(\Div(\Delta) \cap N) \preceq_L \vee_{\!\scriptscriptstyle L}(\Div(\Delta) \cap N)$, thus $\vee_{\!\scriptscriptstyle L}(\Div(\Delta) \cap N) = \vee_{\!\scriptscriptstyle R}(\Div(\Delta) \cap N)$, because $M$ is atomic. 
We denote by $\Delta_N$ that element. 
By definition, we have $\Div(\Delta) \cap N \subseteq \Div(\Delta_N)$. 
On the other hand, we have $\Div(\Delta_N) \subseteq \Div(\Delta) \cap N$, since $\Div(\Delta) \cap N$ contains $\Delta_N$. 
So, $\Div(\Delta_N) = \Div(\Delta) \cap N$, and $\Delta_N$ is balanced. 
It remains to show that $\Div(\Delta_N)$ generates $N$. 
Let $a \in N$. 
Let $a=a_1 a_2 \cdots a_n$ be its left greedy normal form. 
Then $a_i \in \Div(\Delta)$ by definition, and $a_i \in N$ since $N$ is special, thus $a_i \in \Div(\Delta) \cap N = \Div(\Delta_N)$. 
The second part of the lemma is left to the reader.
\qed

\bigskip\noindent
{\bf Proposition 2.5.}
{\it Any Garside monoid $M$ satisfies the following properties.
\begin{itemize}
\item[(P1)] 
The canonical morphism~$\iota: M\to G(M)$ is into.
\item[(P2)] 
If $N$ is a parabolic submonoid, then the associated standard parabolic subgroup is isomorphic to $G(N)$. Moreover, one has $G(N)\cap M = N$ in $G(M)$.
\item[(P3)] 
If $N$ and $N'$ are parabolic submonoids, then $N\cap N'$ is parabolic and $G(N\cap N') = G(N)\cap G(N')$.
\item[(P4)] 
The group $G(M)$ is torsion free.
\end{itemize}}

\bigskip\noindent
{\bf Proof.}
Property (P1) is proved in \cite{DePa}. Property (P2) is implicit in Lemma 2.4 and proved in \cite{God}. Property (P3) is also proved in \cite{God}. Property (P4) is proved in \cite{Deh2}.
\qed

\subsection{Presentations}

\bigskip\noindent
{\bf Definition.}
Recall that an (undirected simple) \emph{graph} is an ordered pair $\Gamma = (S(\Gamma), E(\Gamma)) = (S,E)$ consisting of a set $S$ of \emph{vertices} together with a set $E$ of \emph{edges}, that are $2$-element subsets of $S$.
With $\Gamma$ we associate the set
\[
\hat E(\Gamma) = \{(a,b) \in S \times S \mid \{a,b\} \in E(\Gamma) \}\,.
\]
A \emph{partial complement} on $S$ (based on the graph $\Gamma$) is a mapping $f: \hat E(\Gamma) \to S^*$. The \emph{monoid associated to $f$ on the left} is the monoid $M_L(\Gamma,f)$ defined by the following monoid presentation
\[
M_L(\Gamma,f)=\langle S \mid x\,f(x,y) = y\,f(y,x) \text{ for } \{x,y\} \in E(\Gamma) \rangle^+\,.
\]
Similarly, the  \emph{monoid associated to $f$ on the right} is the monoid $M_R(\Gamma,f)$ defined by the monoid presentation
\[
M_R(\Gamma,f)=\langle S \mid f(y,x)\,x = f(x,y)\,y \text{ for } \{x,y\} \in E(\Gamma) \rangle^+\,.
\]

\bigskip\noindent
{\bf Definition.}
Let $M$ be a preGarside monoid, and let $S$ be a generating set for $M$ which does not contain $1$. Let $\Gamma_L = \Gamma_L(S,M)$ denote the graph on $S$ such that $\{a,b\} \in E(\Gamma_L)$ if and only if $a \neq b$ and $a \vee_{\!\scriptscriptstyle L} b$ exists. Similarly, $\Gamma_R = \Gamma_R(S,M)$ denotes the graph on $S$ such that $\{a,b\} \in E(\Gamma_R)$ if and only if $a \neq b$ and $a \vee_{\!\scriptscriptstyle R} b$ exists. A \emph{left selector} on $S$ in $M$ is a partial complement $f_L$ on $S$ based on $\Gamma_L$ such that $x\,f_L(x,y)$ and $y\,f_L(y,x)$ represent $x\vee_{\!\scriptscriptstyle L} y$ for all $\{x,y\} \in E(\Gamma_L)$. Similarly, a \emph{right selector} on $S$ in $M$ is a partial complement $f_R$ on $S$ based on $\Gamma_R$ such that $f_R(y,x)\,x$ and $f_R(x,y)\,y$ represent $x\vee_{\!\scriptscriptstyle R} y$ for all $\{x,y\} \in E(\Gamma_R)$.

\bigskip\noindent
The following theorem extends \cite[Thm.4.1]{DePa} and is proved in the same way.

\bigskip\noindent
{\bf Theorem 2.6.}
{\it Let $M$ be a preGarside monoid, let $S$ be a generating set for $M$ that does not contain $1$, let $\Gamma_L=\Gamma_L(S,M)$ be as defined above, and let $f_L$ be a left selector on $S$ in $M$. Then $M\simeq M_L(\Gamma_L,f_L)$.}

\bigskip\noindent
{\bf Proof.}
We denote by $\equiv$ the congruence relation on $S^*$ such that $M = (S^*/ \equiv)$. On the other hand, we denote by $\equiv_L$ the congruence relation on $S^*$ generated by the pairs $(x\,f_L(x,y), 
\allowbreak
y\,f_L(y,x))$, $(x,y) \in \hat E(\Gamma_L)$. If $(x,y) \in \hat E(\Gamma_L)$, then $x\,f_L(x,y)$ and $y\, f_L(y,x)$ represent the same element, $x \vee_{\!\scriptscriptstyle L} y$, thus $x\,f_L(x,y) \equiv y\,f_L(y,x)$. So, if $u \equiv_L v$, then $u \equiv v$, for all $u,v \in S^*$.

\bigskip\noindent
For $w \in S^*$, we denote by $\overline{w}$ the element of $M$ represented by $w$. Let $\nu: M \to \N$ be a norm. We take $u,v \in S^*$ such that $u \equiv v$, and turn to prove by induction on $\nu(\overline{u}) = \nu(\overline{v})$ that $u \equiv_L v$. Set $\epsilon =()$, the empty word. If $\nu(\overline{u}) = \nu(\overline{v})=0$, then $u=v=\epsilon$, thus $u \equiv_L v$. Suppose that $\nu(\overline{u}) = \nu(\overline{v})>0$ plus the induction hypothesis. Write $u=xu'$ and $v=yv'$, where $x,y \in S$ and $u',v' \in S^*$. If $x=y$, then, by cancellativity, $u' \equiv v'$, thus, by the induction hypothesis, $u' \equiv_L v'$, therefore $u=xu' \equiv_L v=xv'$. Hence, we can suppose $x \neq y$. Since $x,y \preceq_L \overline{u}$, the element $x \vee_{\!\scriptscriptstyle L} y$ exists, and $(x \vee_{\!\scriptscriptstyle L} y) \preceq_L \overline{u}$. Choose $w \in S^*$ such that $\overline{u} = (x \vee_{\!\scriptscriptstyle L} y)\overline{w}$. By cancellativity, we have $f_L(x,y) w \equiv u'$, thus, by the induction hypothesis, $f_L(x,y)w \equiv_L u'$, therefore $xf_L(x,y)w \equiv_L u$. Similarly, $yf_L(y,x)w \equiv_L v$. Since $xf_L(x,y)w \equiv_L yf_L(y,x)w$, we conclude that $u \equiv_L v$.
\qed

\bigskip\noindent
{\bf Remark.}
In all the algorithmic studies in the theory of Garside groups, a Garside monoid (or group) is given by a finite generating set $S$ together with two complements $f_L,f_R$ on $S$ such that $M=M_L(K_S,f_L)=M_R(K_S,f_R)$, where $K_S$ denotes the complete graph on $S$. There is no algorithm that, given a finite set $S$ and two complements $f_L,f_R$ on $S$ (based on $K_S$) such that $M_L(K_S,f_L)=M_R(K_S,f_R)$, determines whether $M_L(K_S,f_L)$ is a Garside monoid. However, there are partial algorithms, say methods, to solve this question (see \cite{Deh}, for instance). Anyway, it seems reasonable to us that, in order to study algorithmic questions on preGarside monoids and groups, one has to start with a finite set $S$ and two complements $f_L,f_R$ based on $\Gamma_L,\Gamma_R$, respectively, and to assume that $M_L(\Gamma_L,f_L) = M_R(\Gamma_R,f_R)$ is a preGarside monoid.

\bigskip\noindent
Recall that the set of atoms of an atomic monoid $M$ is denoted by $\AA(M)$. It is easily seen that, if $M$ is a preGarside monoid and $N$ is a parabolic submonoid, then $\AA(N) \subset \AA(M)$. The proof of the following is left to the reader.

\bigskip\noindent
{\bf Lemma 2.7.}
{\it Let $M$ be a preGarside monoid. Let $S$ be a generating set for $M$, let $\Gamma_L = \Gamma_L(M,S)$, and let $f_L$ be a left-selector on $S$ in $M$. An element $x \in S$ is an atom if and only if, for all $y \in S \setminus \{x\}$, either $\{x,y\} \not\in E(\Gamma_L)$, or $f_L(x,y) \neq \epsilon$.}
\qed

\bigskip\noindent
So, without loss of generality, in order to study algorithmic questions on preGarside groups, one may assume that the generating set $S$ of $M=M(\Gamma_L,F_L)=M(\Gamma_R,f_R)$ is the set of atoms. Now, the following shows that (if $S$ is finite) there is an effective way to determine all parabolic submonoids of a preGarside monoid.

\bigskip\noindent
{\bf Theorem 2.8.}
{\it Let $M$ be a preGarside monoid, let $f_L$ be a left selector on $\AA(M)$ in $M$, and let $f_R$ be a right selector on $\AA(M)$ in $M$. Let $X$ be a subset of $\AA(M)$, and let $N$ be the submonoid of $M$ generated by $X$. Then $N$ is a parabolic submonoid if and only if the following properties hold. 
\begin{itemize}
\item[(a)]
For all $x,y \in X$, $x \neq y$, if $x \vee_{\!\scriptscriptstyle L} y$ exists, then $f_L(x,y), f_L(y,x) \in X^*$. Similarly, for all $x,y \in X$, $x \neq y$, if $x \vee_{\!\scriptscriptstyle R} y$ exists, then $f_R(x,y), f_R(y,x) \in X^*$.
\item[(b)]
For all $x \in X$ and $y \in \AA(M) \setminus X$, if $x \vee_{\!\scriptscriptstyle L} y$ exists, then $f_L(x,y) \not\in X^*$. Similarly, for all $x \in X$ and $y \in \AA(M) \setminus X$, if $x \vee_{\!\scriptscriptstyle R} y$ exists, then $f_R(y,x) \not\in X^*$.
\end{itemize}}

\bigskip\noindent
{\bf Proof.}
Assume that $N$ is parabolic. 
Let $x,y \in X$, $x \neq y$. 
If $x \vee_{\!\scriptscriptstyle L} y$ exists, then $x \vee_{\!\scriptscriptstyle L} y \in N$, and any expression of $x \vee_{\!\scriptscriptstyle L} y$ belongs to $X^*$, thus $f_L(x,y), f_L(y,x) \in X^*$. 
Similarly, $f_R(x,y), f_R(y,x) \in X^*$ if $x \vee_{\!\scriptscriptstyle R} y$ exists. 
Let $x \in X$ and $y \in \AA(M) \setminus X$. Suppose that $x \vee_{\!\scriptscriptstyle L} y$ exists and that $f_L(x,y) \in X^*$. Then $x \vee_{\!\scriptscriptstyle L} y = \overline{x\,f_L(x,y)} \in N$, $y \preceq_L x \vee_{\!\scriptscriptstyle L} y$, and $y \not\in N$: a contradiction. So, $f_L(x,y) \not\in X^*$. Similarly, $f_R(y,x) \not\in X^*$ if $x \vee_{\!\scriptscriptstyle R} y$ exists.

\bigskip\noindent
Now, we assume that $N$ satisfies Properties (a) and (b), and turn to prove that $N$ is parabolic. Firstly, we take $a \in N$ and an expression $a=y_1 \cdots y_m$, with $y_j \in \AA(M)$ for all $j \in \{1, \dots, m\}$, and we prove that $y_j \in X$ for all $j \in \{1, \dots, m\}$. We take a norm $\nu: M \to \N$ and we argue by induction on $\nu(a)$. The case $\nu(a)=0$ being trivial, we can assume that $\nu(a) >0$ plus the induction hypothesis.

\bigskip\noindent
We choose an expression $a=x_1 \cdots x_n$ of $a$ such that $x_i \in X$ for all $i \in \{1, \dots, n\}$. 
Since $x_1 \preceq_L a$ and $y_1 \preceq_L a$, $x_1 \vee_{\!\scriptscriptstyle L} y_1$ exists. 
Suppose first that $x_1=y_1$. 
Then $\overline{x_2 \cdots x_n} = \overline{y_2 \cdots y_m}$, thus, by the induction hypothesis, $y_2, \dots, y_m$ belong to $X$. 
Now, suppose that $y_1 \neq x_1$. 
There exist $z_1, \dots, z_p$ in $\AA(M)$ such that
\[
a = \overline{x_1\, f_L(x_1,y_1)\, z_1 \cdots z_p} = \overline{y_1\, f_L(y_1,x_1)\,z_1 \cdots z_p}\,.
\]
We have
\[
\overline{f_L(x_1,y_1)\,z_1 \cdots z_p} = \overline{x_2 \cdots x_n} \in N\,,
\]
thus, by the induction hypothesis, $f_L(x_1,y_1)$ belongs to $X^*$ and $z_1, \dots, z_p$ lie in $X$. 
In particular, by Property~(b), we have $y_1 \in X$, and so, by Property (a), $f_L(y_1,x_1)$ belongs to $X^*$.  
Finally,
\[
\overline{y_2 \cdots y_m} = \overline{f_L(y_1,x_1)\,z_1 \cdots z_p} \in N,
\]
thus, by the induction hypothesis, $y_2, \dots, y_m$ lie in $X$.

\bigskip\noindent
Now, we take $a,b \in N$ such that $a \vee_{\!\scriptscriptstyle L} b$ exists, and turn to show that $a \vee_{\!\scriptscriptstyle L} b \in N$. We argue by induction on $\nu(a \vee_{\!\scriptscriptstyle L} b)$. The case $\nu(a \vee_{\!\scriptscriptstyle L} b) =0$ being trivial, we may assume that $\nu(a \vee_{\!\scriptscriptstyle L} b) >0$ plus the induction hypothesis. If $a=1$, then $a \vee_{\!\scriptscriptstyle L} b = b \in N$. Similarly, if $b=1$, then $a \vee_{\!\scriptscriptstyle L} b = a \in N$. So, we can assume that $a \neq 1$ and $b \neq 1$. Write $a=xa_1$ and $b=yb_1$, where $x,y \in X$ and $a_1, b_1 \in N$. We have $xa_1 = a \preceq_L (a \vee_{\!\scriptscriptstyle L} b)$ and $\overline{x\,f_L(x,y)} = (x \vee_{\!\scriptscriptstyle L} y) \preceq_L (a \vee_{\!\scriptscriptstyle L} b)$, thus $a_1 \vee_{\!\scriptscriptstyle L} \overline{f_L(x,y)}$ exists and $x(a_1 \vee_{\!\scriptscriptstyle L} \overline{f_L(x,y)}) \preceq_L (a \vee_{\!\scriptscriptstyle L} b)$. By the induction hypothesis, it follows that $(a_1 \vee_{\!\scriptscriptstyle L} \overline{f_L(x,y)}) \in N$. Similarly, $b_1 \vee_{\!\scriptscriptstyle L} \overline{f_L(y,x)}$ exists, $y(b_1 \vee_{\!\scriptscriptstyle L} \overline{f_L(y,x)}) \preceq_L (a \vee_{\!\scriptscriptstyle L} b)$, and $(b_1 \vee_{\!\scriptscriptstyle L} \overline{f_L(y,x)}) \in N$. Since, by the above, any expression of $a_1 \vee_{\!\scriptscriptstyle L} \overline{f_L(x,y)}$ lies in $X^*$, there is $a_2 \in N$ such that $\overline{f_L(x,y)}a_2 = (a_1 \vee_{\!\scriptscriptstyle L} \overline{f_L(x,y)})$. Similarly, there is $b_2 \in N$ such that $\overline{f_L(y,x)}b_2 = (b_1 \vee_{\!\scriptscriptstyle L} \overline{f_L(y,x)})$. We have 
\[
(x \vee_{\!\scriptscriptstyle L} y)a_2 = \overline{x\,f_L(x,y)}a_2 = x (a_1 \vee_{\!\scriptscriptstyle L} \overline{f_L(x,y)}) \preceq_L (a \vee_{\!\scriptscriptstyle L} b)\,.
\]
Similarly, $(x \vee_{\!\scriptscriptstyle L} y)b_2 \preceq_L (a \vee_{\!\scriptscriptstyle L} b)$, thus $a_2 \vee_{\!\scriptscriptstyle L} b_2$ exists, and $(x \vee_{\!\scriptscriptstyle L} y)(a_2 \vee_{\!\scriptscriptstyle L} b_2) \preceq_L (a \vee_{\!\scriptscriptstyle L} b)$. By the induction hypothesis, it follows that $(a_2 \vee_{\!\scriptscriptstyle L} b_2) \in N$. Finally, 
\begin{gather*}
a = xa_1 \preceq_L x(a_1 \vee_{\!\scriptscriptstyle L} \overline{f_L(x,y)}) = (x \vee_{\!\scriptscriptstyle L} y) a_2 \preceq_L (x\vee_{\!\scriptscriptstyle L} y)(a_2 \vee_{\!\scriptscriptstyle L} b_2)\,,\\
b = yb_1 \preceq_L y(b_1 \vee_{\!\scriptscriptstyle L} \overline{f_L(y,x)}) = (x \vee_{\!\scriptscriptstyle L} y) b_2 \preceq_L (x\vee_{\!\scriptscriptstyle L} y)(a_2 \vee_{\!\scriptscriptstyle L} b_2)\,,
\end{gather*}
thus $(a \vee_{\!\scriptscriptstyle L} b) \preceq (x \vee_{\!\scriptscriptstyle L} y)(a_2 \vee_{\!\scriptscriptstyle L} b_2)$, therefore
\[
(a \vee_{\!\scriptscriptstyle L} b) = (x \vee_{\!\scriptscriptstyle L} y)(a_2 \vee_{\!\scriptscriptstyle L} b_2) \in N\,.
\]
It is easily proved in the same way that $a \vee_{\!\scriptscriptstyle R} b$ lies in $N$ if it exists.
\qed



\section{Amalgamation of monoids}

Given two groups $G_1$ and $G_2$ with a common subgroup $H$, it is known since the work of Schreier \cite{Sch} (see also \cite{Eps}) that both groups $G_1$ and $G_2$ embed in the amalgamated product $G_1\!*_H\! G_2$ above $H$. In this context, given transversals of $G_1/H$ and $G_2/H$ that contain $1$, it is also known by results of Serre \cite{Ser} that every element of the amalgamated product $G_1\!*_H\! G_2$ has a unique amalgam decomposition. In the context of monoids this is not true anymore (see \cite{How,How2}). In particular, amalgam decompositions do not exist in general, and, even if we can effectively decide whether an element of $M_i$ lies in $N$ for $i=1,2$, the word problem may be not decidable in $M_1\!*_N\! M_2$ (see \cite{BiMaMe}). The aim of this section is to prove that, if $M_1$ and $M_2$ are preGarside monoids and $N$ is a common parabolic submonoid, then $M_1\!*_N\! M_2$ is also a preGarside monoid, $M_1$ and $M_2$ embed in $M_1\!*_N\! M_2$, and amalgam decompositions exist in the later monoid. In Section 4 we will use this to define and investigate the notion of preGarside monoids and groups of FC type.

\subsection{Special amalgamation of monoids}

If $M_1$, $M_2$ are two monoids, we denote by $F^+(M_1,M_2)$ the semigroup $\{(g_1,\dots,g_n)\mid n\geq 1,\ g_i\in M_1\cup M_2\}$ equipped with the concatenation operation. For $g=(g_1,\dots,g_n)$ we set $|g| = n$, and we define $\varepsilon_i(g)$ by $g_i\in M_{\varepsilon_i(g)}$ for $1\leq i\leq |g|$. Note that, with this definition, the semigroup $F^+(M_1,M_2)$ is not a monoid, since it does not contain the empty sequence $\epsilon=()$. 

\bigskip\noindent
{\bf Definition.}
Let $M_1, M_2$ and $N$ be three monoids such that there exist injective morphisms of monoids $\iota_1:N\to M_1$ and $\iota_2:N\to M_2$. The \emph{amalgamated  product of the monoids $M_1$ and $M_2$ over $N$} is the monoid $M_1\!*_N\!M_2$ obtained as the quotient of the free semigroup $F^+(M_1,M_2)$ by the congruence $\equiv$ generated by the binary relation $\equiv_0$ defined by 
\[
(g_1,\dots,g_n)\equiv_0 (g_1,\dots,g_{i-1},\tilde g_i,g_{i+2},\dots,g_n)
\]
if one of the following conditions holds: 
\begin{itemize}
\item[(a)]
$\tilde g_i  = g_ig_{i+1}$, with $\varepsilon_i = \varepsilon_{i+1}$;
\item[(b)]
$\tilde g_i  = g_i\,(\iota_{\varepsilon_i}\!\circ\!\iota_{\varepsilon_{i+1}}^{-1})(g_{i+1})$, with $\varepsilon_i \neq \varepsilon_{i+1}$ and $g_{i+1}$ belonging to $\iota_{\varepsilon_{i+1}}(N)$;
\item[(c)]
$\tilde{g}_i= (\iota_{\varepsilon_{i+1}}\!\circ\!\iota_{\varepsilon_i}^{-1})(g_i)\,g_{i+1}$, with $\varepsilon_i \neq \varepsilon_{i+1}$ and $g_i$ belonging to $\iota_{\varepsilon_i}(N)$;
\end{itemize}
where $\varepsilon_i = \varepsilon_i(g_1,\dots,g_{n})$ and  $\varepsilon_{i+1} = \varepsilon_{i+1}(g_1,\dots,g_{n})$.

\bigskip\noindent
{\bf Definition.}
We say that $M_1\!*_N\!M_2$ is a \emph{special amalgam} when $M_1, M_2$ and $N$ are three monoids with two injective morphisms of monoids $\iota_1:N\to M_1$ and $\iota_2:N\to M_2$ such that $\iota_1(N)$ and  $\iota_2(N)$ are special submonoids of $M_1$ and $M_2$, respectively.

\bigskip\noindent
{\bf Proposition 3.1.}
{\it Let $M_1\!*_N\!M_2$ be a special amalgam.
Then the canonical morphisms $j_1: M_1 \to M_1\!*_N\!M_2$ and  $j_2: M_2\to M_1\!*_N\!M_2$ are injective, and the submonoids $j_1(M_1)$, $j_2(M_2)$ and $j_1\circ\iota_1(N)$  are special in $M_1\!*_N\!M_2$. 
Moreover, $j_1(M_1)\cap j_2(M_2) = (j_1\circ\iota_1)(N) = (j_2\circ\iota_2)(N)$.}

\bigskip\noindent
{\bf Proof.}
Let $g$ belong to $M_1$ and assume $(g)\equiv (g_1,\dots,g_n)$. Using that $\iota_1(N)$ and $\iota_2(N)$ are special, we prove by an easy induction on the number of elementary relations $\equiv_0$ needed to transform $(g)$ into $(g_1,\dots,g_n)$ that, firstly, for each $i$, $g_i$ belongs either to $M_1$, or to $\iota_2(N)$ and, secondly, $g = h_1\cdots h_n$ in $M_1$, where $h_i = g_i$,  or $h_i =\iota_1\circ\iota_2^{-1}(g_i)$. Therefore, the morphism $j_1$ is injective, and $j_1(M_1)$ is a special submonoid. The rest of the proposition follows from similar arguments.
\qed

\bigskip\noindent
Note that the fact that the canonical morphisms $j_1: M_1 \to M_1\!*_N\!M_2$ and  $j_2: M_2\to M_1\!*_N\!M_2$ are injective was known before \cite{How}, because a special submonoid is \emph{unitary} (see~\cite[p.~273]{BiMaMe} for a definition).  
In the sequel, when $M_1\!*_N\!M_2$ is a special amalgam, we identify $M_1, M_2$ and $N$ with their images in $M_1\!*_N\!M_2$.

\subsection{Amalgam decomposition}

In order to introduce the second main notion of this section, we need first to recall the notion of a \emph{confluent reduction rule}. Consider a set $X$. A \emph{reduction rule} on $X$ is a map $f$ from a set $Y$ to the set $\PP (X\times X)$ of subsets of $X\times X$ such that, for every $y$ in $Y$ and every $x$ in $X$, there is at most one $x'$ in $X$ such that $(x,x')$ belongs to $f(y)$. In this case, we write $x' \stackrel{y}{\longleftarrow} x$, and we say that \emph{$x$ reduces to $x'$ by a reduction of type $y$}. We denote by  $\stackrel{*}{\longleftarrow}$ the reflexive-transitive binary relation induced by $f$. In other words, $x' \stackrel{*}{\longleftarrow} x$ if there is a finite sequence of reductions
\[
x' = x_n \stackrel{y_n}{\longleftarrow}\cdots \stackrel{y_2}{\longleftarrow}x_1\stackrel{y_1}{\longleftarrow} x_0 = x\,.
\]
Finally, by $\longleftrightarrow$ we denote the induced equivalence relation on $X$. We say that the reduction rule is \emph{globally confluent} if, for every two elements $x_1,x_2$ in $X$ such that $x_1\longleftrightarrow x_2$, there exists $z$ in $X$ such that $ z\stackrel{*}{\longleftarrow}x_1$ and $z\stackrel{*}{\longleftarrow} x_2$. We say that the reduction rule is \emph{locally confluent} if the following property holds: for all $x_1,x_2, x_3$ in $X$ such that $x_1\stackrel{y_1}{\longleftarrow}x_3$ and $x_2\stackrel{y_2}{\longleftarrow} x_3$, there exists $x_4$ in $X$ such that $x_4 \stackrel{*}{\longleftarrow} x_1$ and $x_4 \stackrel{*}{\longleftarrow} x_2$.  Finally, we say that the the reduction rule is \emph{Noetherian} if every infinite sequence $(x_i)$ such that $ x_{i+1} \stackrel{*}{\longleftarrow} x_i$ has to stabilize. The following is classical in the subject.

\bigskip\noindent
{\bf Lemma 3.2}
(Diamond Lemma).
{\it Every Noetherian and locally confluent reduction rule is globally confluent.}
\qed

\bigskip\noindent
Note that a consequence of the Diamond lemma is that, for a Noetherian and locally confluent reduction rule $\longleftarrow$, every equivalence class for the relation $\longleftrightarrow$ possesses a unique minimal element. In other words, if $\CC$ is an equivalence class for $\longleftrightarrow$, there is a unique $x$ in $\CC$ such that $x \stackrel{*}{\longleftarrow} y$ for all $y$ in $\CC$.

\bigskip\noindent
Now, recall some classical notions from semigroup theory \cite{Gre}. Let $M$ be a monoid, and let $N$ be a submonoid. We define two relations $\RR_N$ and $\LL_N$ on $M$ setting $g\RR_N h$ and $g\LL_N h$ if $gN = hN$ and $Ng = Nh$, respectively. When $N = M$, the relations $\RR_M$ and $\LL_M$ are denoted by $\RR$ and $\LL$, respectively. They are called the \emph{Green relations} on $M$. In what follows, we assume the set $\PP(M)$ of subsets of $M$ to be endowed with the reduction rule $P_1 \longleftarrow P_2$ if $P_1\varsupsetneq P_2$. This induces a reduction rule on the set of $\RR_N$-classes, $\{gN, g\in M\}$, that verifies $g_1N\longleftarrow g_2N$ if $g_2\in g_1N$ and $g_1\notin g_2N$. Similarly, it induces a reduction rule on the set of $\LL_N$-classes, $\{Ng, g\in M\}$.

\bigskip\noindent
{\bf Definition.}
Let $M$ be a monoid, and let $N$ be a submonoid. We say that $N$ has the \emph{$\LL$ confluence property} if the reduction rule $\longleftarrow$ is Noetherian and locally confluent on the set of $\LL_N$-classes. Similarly, We say that $N$ has the \emph{$\RR$ confluence property} if the reduction rule $\longleftarrow$ is Noetherian and locally confluent on the set of $\RR_N$-classes, and we say that $N$ has the \emph{confluence property} if it has both $\LL$ and $\RR$ confluence properties. 

\bigskip\noindent
{\bf Remark.}
If $N$ is a special submonoid of a monoid $M$, then $N$ is minimal for the reduction rule $\longleftarrow$ in the set of $\LL_N$-classes, as well as in the set of $\RR_N$-classes.

\bigskip\noindent
If $M$ is a cancellative monoid, $N$ is a special submonoid with the confluence property, and $T$ is a set of representatives of the minimal $\RR_N$-classes that contains $1$, then, for every element $g$ of $M$, there exists a unique pair $(g_1,h)$ in $T\times N$ such that $g = g_1h$. In the sequel we will set $[g]_N = (g_1,h)$ if $g_1\neq 1$ and $[g]_N = (h)$ otherwise.

\bigskip\noindent
{\bf Proposition 3.3.}
{\it Let $M_1\!*_N\!M_2$ be a special amalgam such that $N$ has the confluence property in both, $M_1$ and $M_2$. For $i = 1,2$, consider a set $T_i$ of representatives of the minimal $\RR_N$-classes in $M_i$ that contains $1$. If $M_1$ and $M_2$ are cancellative, then every element $g$ of $M_1\! *_N\!M_2$ has a unique decomposition 
\[
g = g_1\cdots g_m h
\]
such that each $g_i$ belongs to either $T_1\setminus\{1\}$ or $T_2\setminus\{1\}$, $h$ belongs to $N$, and two consecutive $g_i$ do not lie in the same set of representatives.}

\bigskip\noindent
In the sequel the sequence $(g_1,\dots,g_m,h)$ is called the \emph{(left) amalgam decomposition} of $g$ and is denoted by $[g]_N$. We denote the integer $m$ by $\ell_N(g)$. If $g$ belongs to $N$, we set $\ell_N(g) = 0$. Clearly, $\ell_N(g)$ does not depend on the choice of the transversals. Note that, choosing a set $T_i$ of representatives of the minimal $\LL_N$-classes in $M_i$ that contains $1$, one can define similarly a \emph{right amalgam decomposition} $(h,g_m,\dots, g_1)$ of every $g$. Moreover, we have $m =\ell_N(g)$. Proposition~3.3 is a consequence of the following lemma.

\bigskip\noindent
{\bf Lemme 3.4.}
{\it Let $M_1\!*_N\!M_2$ be a special amalgam such that $N$ has the confluence property in both, $M_1$ and $M_2$. For $i = 1,2$, we choose a set $T_i$ of representatives of the minimal $\RR_N$-classes in $M_i$ that contains $1$. Assume $M_1$ and $M_2$ are cancellative. Consider the reduction rule on $F^+(M_1,M_2)$ whose types are in $\{a,b,c\}\times \N^*$ and which is defined in the following way. If $X= (g_1,\dots, g_m)$ and $X' = (g'_1,\dots, g'_n)$, then
\begin{itemize}
\item[(a)]
$X'\stackrel{(a,i)}{\longleftarrow} X$ if $m = n$, $1\leq i\leq m-1$, $(g_i,g_{i+1})\in(M_1\!\times\! M_1)\cup(M_2\!\times\! M_2)$, $(g'_i,g'_{i+1}) = [g_ig_{i+1}]_N$, and $g'_j = g_j$ for $j\neq i,i+1$;
\item[(b)]
$X'\stackrel{(b,i)}{\longleftarrow}X$ if $n= m+1$, $1\leq i\leq m$, $g_i\not\in N\cup T_1\cup T_2$, $(g'_i,g'_{i+1}) = [g_i]_N$, and $g'_j = g_j$, $g'_{k+1} = g_{k}$ for $j< i < k$;
\item[(c)]
$X'\stackrel{(c,i)}{\longleftarrow}X$ if $n= m-1$, $1\leq i\leq m-1$, $(g_i,g_{i+1})\in N\times N$, $g'_i = g_ig_{i+1}$, and $g'_j = g_j$, $g'_k = g_{k+1}$ for $j< i < k$.
\end{itemize}
Then the reduction rule is Noetherian and locally confluent.}

\bigskip\noindent
{\bf Proof.}
Consider a sequence $(X_k)_k$ in $F^+(M_1,M_2)$ such that $X_{k+1}\stackrel{*}{\longleftarrow} X_k$ for all $k$. If $X = (g_1,\dots, g_n)$ belongs to $F^+(M_1,M_2)$, set
\begin{gather*}
|X|_{\not\in N} = |\{i\mid g_i\not\in N\}|\,,\quad
|X|_{\not\in T\cup N} = |\{i\mid g_i\not\in T_1\cup T_2\cup N\}|\,,\\
\Inv(X) = |\{(i,j)\mid i<j,\ g_i\in N,\ g_j \not\in N \}|\,.
\end{gather*}
It is easily seen that $|X_k|_{\not\in N} \ge |X_{k+1}|_{\not\in N}$ and $|X_k|_{\not\in T \cup N} \ge |X_{k+1}|_{\not\in T \cup N}$ for all $k \in \N$, therefore there exists $K \in \N$ such that the sequences $|X_k|_{\not\in N}$ and $|X_k|_{\not\in T \cup N}$ stabilize for $k\geq K$. It follows that, for $k\geq K$, the only reduction rules that can be applied are either of type $(c,i)$ or of type $(a,i)$. Moreover, in the latter case, $(g_i,g_{i+1})$ has to belong either to $N\times T_{\varepsilon}$ or to $T_{\varepsilon}\times N$, where $\varepsilon = 1,2$. Now, $|X_k| \ge |X_{k+1}|$ and $\Inv(X_k) \ge \Inv(X_{k+1})$ for $k\geq K$, thus there exists $K_1\geq K$ such that the sequences $|X_k|$ and $\Inv(X_k)$ stabilize for $k\geq K_1$. For $k\geq K_1$, the only reduction rules that can be applied are of type $(a,i)$ with $(g_i,g_{i+1})$ in $T_{\varepsilon}\times N$. But, in this case, $[g_ig_{i+1}]_N = (g_i,g_{i+1})$. Therefore, $(X_k)_k$ stabilizes for $k\geq K_1$. This shows that the reduction rule is Noetherian.

\bigskip\noindent
Now, assume $X' \stackrel{(e,i)}{\longleftarrow} X$ and $X''\stackrel{(f,j)}{\leftarrow} X$ with $X = (g_1,\dots, g_m)$ and $1\leq i \leq j\leq m$. We need to find some $X'''$ such that $X''' \stackrel{*}{\longleftarrow} X'$ and $X''' \stackrel{*}{\longleftarrow} X''$. If either $j\geq i+2$, or $e = b$ with $j\geq i+1$, then such a $X'''$ is easily found. Also, if $e=f=c$ and $j=i+1$, then $X'''$ is easily constructed. So, we should treat the following remaining three cases: (1) $\{e,f\} = \{a,b\}$ and $i = j$; (2) $e = c$, $f = a$ and $j = i+1$; (3) $e = a$, $j = i+1$ and $f= a,b,c$.

\bigskip\noindent
{\it Case 1.}
Assume by symmetry that $e = a$ and $f = b$. Then we have 
\[
X' \stackrel{(c,i+1)}{\longleftarrow}\cdot\stackrel{(a,i)}{\longleftarrow}\cdot \stackrel{(a,i+1)}{\longleftarrow} X''\,
\]
by the uniqueness of the decomposition $[g_ig_{i+1}]_N$.

\bigskip\noindent
{\it Case 2.}
We have 
\[
X''' \stackrel{(a,i)}{\longleftarrow} X' \quad \text{and} \quad X'''  \stackrel{(c,i+1)}{\longleftarrow}\cdot\stackrel{(a,i)}{\longleftarrow} X''
\]
by the uniqueness of the decomposition $[g_ig_{i+1}g_{i+2}]_N$.

\bigskip\noindent
{\it Case 3.}
Assume first $f = a$. We have two cases depending on whether $g_i$ and $g_{i+2}$ belong to the same $M_j$ or not. In the first case, we have
\[
X'''\stackrel{(c,i+1)}{\longleftarrow}\cdot\stackrel{(a,i)}{\longleftarrow}\cdot \stackrel{(a,i+1)}{\longleftarrow} X'\quad \text{and} \quad X''' \stackrel{(c,i+1)}{\longleftarrow}\cdot\stackrel{(a,i)}{\longleftarrow} X''
\]
by the uniqueness of the decomposition $[g_ig_{i+1}g_{i+2}]_N$. In the second case, $g_{i+1}$ has to lie in $N$. This is the only non-obvious case. We may assume without loss of generality that $m = 2$ and $i = 1$. Up to symmetry, we may also assume $g_1\in M_1$ and $g_3\in M_2$. So, we have $X = (g_1,g_2,g_3)$. Set $X' = (g'_1,g'_2,g'_3)$ and $X'' = (g''_1,g''_2,g''_3)$. Let $X''' = (g'_1,g^{(3)}_2,g^{(3)}_3)$ be such that $X''' \stackrel{(a,2)}{\longleftarrow} X'$. The following equalities have to hold: $g'_3 = g_3$, $[g_1g_2]_N = (g'_1,g'_2)$, and $[g'_2g_3]_N = (g^{(3)}_2,g^{(3)}_3)$. Consider the sequence of reductions
\[
X^{(6)}\stackrel{(c,3)}{\longleftarrow}X^{(5)}\stackrel{(a,2)}{\longleftarrow}X^{(4)}\stackrel{(b,1)}{\longleftarrow} X''\,.
\]
Write $X^{(4)} = (g^{(4)}_1,g^{(4)}_2,g''_2,g''_3)$, $X^{(5)} = (g^{(4)}_1,g^{(5)}_2,g^{(5)}_3,g''_3)$, and $X^{(6)} = (g^{(4)}_1,g^{(5)}_2,g^{(6)}_3)$. We have $[g_1]_N = (g^{(4)}_1,g^{(4)}_2)$. In particular, $g_1g_2 = g'_1g'_2 = g^{(4)}_1g^{(4)}_2g_2$. Thus, $g^{(4)}_1 = g'_1$ and  $g'_2 = g^{(4)}_2g_2$. It follows from the latter equality that, in $M_2$, we have 
\[
g^{(3)}_2g^{(3)}_3 = g'_2g_3 =  g^{(4)}_2g_2g_3 =  g^{(4)}_2g''_2g''_3 = g^{(5)}_2g^{(5)}_3g''_3 = g^{(5)}_2g^{(6)}_3\,.
\]
Hence, $g^{(5)}_2 =  g^{(3)}_2$ and $g^{(3)}_3 =g^{(6)}_3$. In other words, $X^{(6)} = X'''$. Assume now $f = b$. Then 
\[
X' \stackrel{(c,i+1)}{\longleftarrow}\cdot \stackrel{(a,i)}{\longleftarrow} X''
\]
by the uniqueness of the decomposition $[g_ig_{i+1}]_N$. Assume finally $f = c$. Then we have 
\[
X'''\stackrel{(c,i+1)}{\longleftarrow} X' \quad \text{and} \quad X''' \stackrel{(a,i)}{\longleftarrow} X''\,.
\]
\qed

\bigskip\noindent
{\bf Proof of Proposition 3.3.}
We keep the notations of Lemma 3.4. No reduction rule can be applied to $(g_1,\dots, g_n,h)$ if and only if the decomposition $g_1\cdots g_nh$ is as stated in Proposition~3.3. Now, it is immediate from the definitions that the equivalence relation $\longleftrightarrow$ on $F^+(M_1,M_2)$ generated by $\stackrel{*}{\longleftarrow}$ is equal to the relation $\equiv$ given in the definition of a amalgamated product of monoids. Then the result is a consequence of Lemma 3.2 and Lemma 3.4.
\qed

\bigskip\noindent
{\bf Lemma 3.5.}
{\it Let $M_1\!*_N\!M_2$ be a special amalgam such that $M_1$ and $M_2$ are cancellative and $N$ has the confluence property in both, $M_1$ and $M_2$. For $i = 1,2$, consider a set $T_i$ of representatives of the minimal $\RR_N$-classes in $M_i$ which contains $1$. Let $g$ belong to $M_1\!*_N\!M_2$, and assume that $g = g_1\cdots g_n$, where $g_i$ lies in $M_{\varepsilon_i}\setminus N$ and $\varepsilon_i\neq\varepsilon_{i+1}$. Let $g_1', \dots, g_n' \in T_1 \cup T_2$ and $h_0,h_1, \dots, h_n \in N$ be inductively defined by
\[
h_0=1, \quad [h_{i-1}g_i]_N=(g_i',h_i) \text{ if } 1 \le i \le n\,.
\]
Then $(g_1', \dots, g_n',h_n)$ is the amalgam decomposition of $g$.}

\bigskip\noindent
{\bf Proof.}
One has $g = g'_1\cdots g'_n h_n$ and $(g'_1,\dots, g'_n,h_n)$ is an amalgam decomposition: no $g'_i$ is equal to $1$ because $g'_ih_i = h_{i-1}g_i$ and $N$ is special.
\qed

\bigskip\noindent
{\bf Lemma 3.6.}
{\it We keep the notations of Proposition 3.3. Assume $M_1$ and $M_2$ are cancellative. Let $g,g'$ belong to $M_1\!*_{N}\! M_2$, and let $[g]_N = (g_1,\dots, g_n, h)$ and $[g']_N = (g'_1,\dots, g'_m, h')$. If $m\geq 1$, $n\geq 1$, and  $\varepsilon_{m}([g']_N) = \varepsilon_{1}([g]_N)$, then  
\[
[g'g]_N = (g'_1,\dots,g'_{m-1},\tilde{g}_1,\dots,\tilde{g}_n,\tilde{h})\,,
\]
where $(\tilde g_1, \dots, \tilde g_n,\tilde h) = [g_m'h'g]_N$. Otherwise, 
\[
[g'g]_N = (g'_1,\dots,g'_m,\tilde{g}_1,\dots,\tilde{g}_n,\tilde{h})\,,
\]
where $(\tilde g_1, \dots, \tilde g_n, \tilde h) = [h'g]_N$.}

\bigskip\noindent
{\bf Proof.}
This is a direct consequence of Lemma 3.5. In the first case $g'_mhg_1$ does not belong to $N$ because the latter is special and $g'_m$ does not belong to $N$.
\qed

\bigskip\noindent
{\bf Remark.}
Assume $M_1\!*_N\!M_2$ is a special amalgam such that $M_1$ and $M_2$ are cancellative and $N$ has the confluence property in both, $M_1$ and $M_2$. Replacing minimal $\RR_N$-classes by minimal $\LL_N$-classes, and choosing representatives, we can associate to each element $g$ of $M_1\!*_N\!M_2$ a \emph{left amalgam decomposition} $(\tilde{h},\tilde{g}_m,\dots,\tilde{g}_1)$. Then one obtains properties for the left amalgam decompositions that are similar to those proved in  Lemmas 3.5 and 3.6 for right amalgam decompositions.

\bigskip\noindent
{\bf Corollary 3.7.}
{\it Let $M_1\!*_N\!M_2$ be a special amalgam such that $N$ has the confluence property in both $M_1$ and $M_2$. Assume $M_1$ and $M_2$ are cancellative. Then $M_1\!*_N\!M_2$ is cancellative.}

\bigskip\noindent
{\bf Proof.}
Let $g,g',g''$ belong to $M_1\!*_{N}\! M_2$ such that $g'g = g''g$. Clearly, we can assume without loss of generality that $\ell_N(g)\leq 1$. Set $[g']_N = (g'_1,\dots, g'_m, h')$ and  $[g'']_N = (g''_1,\dots, g''_p, h'')$. If $\ell_N(g) = 0$, then $[g'g]_N = (g'_1,\dots, g'_m, h'g)$ and $[g''g]_N = (g''_1,\dots, g''_p, h''g)$. By the uniqueness of the amalgam decomposition, it follows that $m = p$, $g'_1 = g''_1,\dots, g'_p = g''_p$, and $h'g = h''g$ in $N$. The latter equality implies $h' = h''$ because $N$ is cancellative. Now, assume $\ell_N(g) = 1$ (say $g\in M_1\setminus N$). If $g'_m$ lies in $M_2$, then the amalgam decomposition $[g'g]_N$ of $g'g$ is $(g'_1,\dots, g'_{m},\tilde{g}'_{m+1} ,\tilde{h}')$, where $(\tilde{g}'_{m+1} ,\tilde{h}') = [h'g]_N$, and, if $g'_m$ lies in $M_1$, then the amalgam decomposition of $g'g$ is $(g'_1,\dots, g'_{m-1},\tilde{g}'_m ,\tilde{h}')$, where $(\tilde{g}'_m ,\tilde{h}') = [g'_mh'g]_N$. Similarly, if $g''_p$ lies in $M_2$ then the amalgam decomposition $[g''g]_N$ of $g''g$ is $(g''_1,\dots, g''_{p},\tilde{g}''_{p+1} ,\tilde{h}'')$, where $(\tilde{g}''_{p+1} ,\tilde{h}'') = [h''g]_N$, and, if $g''_p$ lies in $M_1$, then the amalgam decomposition of $g''g$ is $(g''_1,\dots, g''_{p-1},\tilde{g}''_p ,\tilde{h}'')$, where $(\tilde{g}''_p ,\tilde{h}'') = [g''_ph''g]_N$. If both $g'_m$ and $g''_p$ belong to $M_2$, then one can easily show that $g' = g''$ in the same way as for the case $\ell_N(g) = 0$. If we assume $g_m'$ lies in $M_1$ and $g''_p$ lies in $M_2$, then the uniqueness of the amalgam decomposition implies that $g'_mh'g = \tilde{g}'_m\tilde{h}' = \tilde{g}''_{p+1}\tilde{h''} = h''g$, thus $g'_mh' = h''$, since $M_1$ is cancellative. This is a contradiction since $N$ is special. Similarly,  we obtain a contradiction if we assume that $g_m'$ lies in $M_2$ and $g''_p$ lies in $M_1$. Assume finally both $g'_m$ and $g''_p$ belong to $M_1$. Again by the uniqueness of the amalgam decomposition, it follows that $m = p$, $g'_1 = g''_1,\dots, g'_{p-1} = g''_{p-1}$, and $g_p'h'g = g''_ph''g$ in $M_1$, which is cancellative. We conclude that $g'_ph' = g''_ph''$ and therefore $g' = g''$.  By similar arguments, the equality $gg' = gg''$ in $M_1\!*_N\!M_2$ implies $g' = g''$.
\qed

\bigskip\noindent
{\bf Corollary 3.8.}
{\it Let $M_1\!*_N\!M_2$ be a special amalgam such that $N$ has the confluence property in both, $M_1$ and $M_2$. Assume $M_1$ and $M_2$ are cancellative and atomic. Then $M_1\!*_N\!M_2$ is atomic.}

\bigskip\noindent
{\bf Proof.}
For all $g$ in $M_i$, $i \in \{1,2\}$, we denote by $\nu_i(g)$ the maximal length of $g$ over the atoms of $M_i$. Recall that the map $\nu_i : M_i \to \N$ is well-defined  and is a norm. Note that, since $N$ is special, we have $\nu_1(h) = \nu_2(h)$ for all $h \in N$. Let $g \in M_1 \! *_N\! M_2$, $g \not\in N$, and let $(g_1, \dots, g_k,h)$ be the amalgam normal form of $g$. Let $\varepsilon(i)$ be the element of $\{ \pm 1\}$ such that $g_i \in M_{\varepsilon(i)}$. We take a sequence $(g_1', \dots, g_k')$ such that $g_i' \in M_{\varepsilon(i)}$ for all $i \in \{1, \dots, k\}$ and $g=g_1' \cdots g_k'$, and turn now to prove that $\nu_{\varepsilon(i)}(g_i')$ is bounded for all $i$. We argue by induction on $k$.

\bigskip\noindent
The case $k=1$ being trivial, we can assume that $k \ge 2$ plus the induction hypothesis. By Lemma 3.5, there exists $h_{k-1} \in N$ such that $h_{k-1} g_k' = g_kh$. This implies that $\nu_{\varepsilon(k)} (g_k')$ is bounded above by $\nu_{\varepsilon(k)}(g_kh)$. Similarly, replacing the right normal form of $g$ by its left normal form, it is shown that $\nu_{\varepsilon(1)}(g_1')$ is bounded. Again by Lemma 3.5, there exists $h_1 \in N$ such that $g_1' = g_1h_1$. Finally, from the equality 
\[
(h_1g_2')g_3' \cdots g_k' = g_2 \cdots g_k h
\]
and the induction hypothesis, it follows that $\nu_{\varepsilon(i)} (g_i')$ is bounded for all $i \in \{3, \dots, k\}$, and that  $\nu_{\varepsilon(2)} (h_1g_2')$ is bounded, too. The fact that $\nu_{\varepsilon(2)} (h_1g_2')$ is bounded implies that $\nu_{\varepsilon(2)} (g_2')$ is bounded. 

\bigskip\noindent
For $g$ in $M_1\!*_N\!M_2$, $g \not\in N$, we set
\[
\nu(g) = \max \{ \nu_{\varepsilon(1)} (g_1') + \cdots + \nu_{\varepsilon(k)}(g_k') \mid g=g_1' \cdots g_k' \text{ et } g_i' \in M_1 \cup M_2 \setminus N \}\,,
\]
and for $g \in N$, we set $\nu(g)=\nu_1(g) = \nu_2(g)$. It is easily checked that $\nu$ is a norm. So, $M_1\!*_N\!M_2$ is atomic.
\qed

\subsection{Amalgamation of preGarside monoids above parabolic submonoids}

The following lemma will allow us to apply the results of the previous subsection to parabolic submonoids in preGarside monoids.

\bigskip\noindent
{\bf Lemma 3.9.}
{\it Let $M$ be a preGarside monoid, and let $N$ be a parabolic submonoid of $M$. Then $N$ has the confluence property. Moreover, each minimal $\RR_N$-class (resp. $\LL_N$-class) has a unique representative.}

\bigskip\noindent
{\bf Proof.}
It is easily checked that the fact that $M$ is atomic and $N$ is special implies that the rewriting rule $\stackrel{*}\longleftarrow$ on the $\RR_N$-classes is Noetherian, and that each minimal class has a unique representative. It remains to show that $\stackrel{*}\longleftarrow$ is locally confluent. Let $g,g',g'' \in M$ such that $g'N \stackrel{*}\longleftarrow gN$ and $g''N \stackrel{*}\longleftarrow gN$. Let $h',h'' \in N$ such that $g=g'h'=g''h''$. We have $h' \preceq_R g$ and $h'' \preceq_R g$, thus $h_0 = h' \vee_{\!\scriptscriptstyle R} h''$ exists and belongs to $N$. Let $g_0 \in M$ such that $g=g_0 h_0$. Using the cancellation property it is then easily shown that $g_0N \stackrel{*}\longleftarrow g'N$ and $g_0N \stackrel{*}\longleftarrow g''N$.
\qed 

\bigskip\noindent
{\bf Remark.}
Let $M_1$ and $M_2$ be preGarside monoids, and let $N$ be a common parabolic submonoid. Then, by Lemma 3.9, left and right amalgam decompositions in $M_1\!*_N\!M_2$ exist, and they are unique in the sense that there is unique choice of transversals $T_1$ and $T_2$ for defining them. 

\bigskip\noindent
{\bf Lemma 3.10.}
{\it Let $M_1$, $M_2$ be preGarside monoids, and let $N$ be a common parabolic submonoid. Let $g,g'$ belong to $M = M_1\!*_N\!M_2$. If there exists $g_0 \in M_1\!*_N\!M_2$ such that $g \preceq_L g_0$ and $g' \preceq_L g_0$, then $g \vee_{\!\scriptscriptstyle L} g'$ exists (in $M_1\!*_N\!M_2$), and $\ell_N(g\vee_{\!\scriptscriptstyle L} g') =  \max(\ell_N(g),\ell_N(g'))$.}

\bigskip\noindent
{\bf Proof.}
We argue by induction on $m =\max(\ell_N(g),\ell_N(g'))$. Set $M=M_1\!*_N\!M_2$. Let $x,x' \in M$ be such that $gx = g'x' = g_0$. Denote by $(t,x_s,\dots,x_1)$ and $(t',x'_{s'},\dots, x'_1)$ the right amalgam decompositions of $x$ and $x'$, respectively. We get $gtx_s\cdots x_1 = g't'x'_{s'}\cdots x'_1=g_0$. Consider the case $m = 0$. Then $g$ and $g'$ belong to $N$. By the uniqueness of the right amalgam decomposition, we get $gt = g't'$ in $N$. Therefore, $g\vee_{\!\scriptscriptstyle L} g'$ exists in $N$. It is easily seen that this element is also the least element in $\{y \in M \mid g \preceq_L y \text{ and } g' \preceq_L y\}$. Consider now the case $m = 1$. The uniqueness of the amalgam decomposition and the existence of a common multiple imply that $g$ and $g'$ both belong either to $M_1$, or to $M_2$. By arguments similar to the previous case, $g\vee_{\!\scriptscriptstyle L} g'$ exists in $M_1$ or $M_2$, and this element is the least element in $\{y \in M \mid g \preceq_L y \text{ and } g' \preceq_L y\}$.

\bigskip\noindent
Now, we assume $m \ge 2$ plus the induction hypothesis. Set $[g]_N = (g_1,\dots,g_k,h)$ and $[g']_N = (g'_1,\dots, g'_{k'},h')$. By Lemma 3.6, the $k-1$ first terms in $[gx]_N$ are $g_1,\dots, g_{k-1}$, and the $k'-1$ first terms in $[g'x]_N$ are $g'_1,\dots, g'_{k'-1}$. Since $gx = g'x'$, it follows that $g_i = g'_i$ for $i\leq \min(k,k')-1$. Hence, upon applying cancellation in the left hand side, we may reduce our study to the case $\min(k,k')\leq 1$. So, we can assume $\ell_N(g)\leq 1$ and $[g']_N = (g'_1,\dots, g'_{k'},h')$ with $k' = m\geq 2$. By the induction hypothesis it follows from the equality $gx = g'_1\cdots g'_mh'x'$  that $g$ and $g'_1\cdots g'_{m-1}$ have a least common multiple $g'' = g \vee_{\!\scriptscriptstyle L} (g'_1\cdots g'_{m-1})$ in $M$ with $\ell_N(g'') = m-1$. Write $x = x_1x_2$ such that $g'' = gx_1 = g'_1\cdots g'_{m-1}\tilde{g}$. One has $\ell_N(\tilde{g})\leq 1$ by Lemma 3.6. By cancellativity, we get $\tilde{g}x_2 = g'_mh'x'$. Applying the case $m = 1$, we deduce that $\tilde{g}$ and $g'_mh'$ have a least common multiple $g''' = \tilde{g} \vee_{\!\scriptscriptstyle L} (g'_mh')$ in $M$ such that $\ell_N(g''') = 1$. Now, $g'_1\cdots g'_{m-1}g'''$ left divides $g_0$ and is the least element in $\{y \in M \mid g \preceq_L y \text{ and } g' \preceq_L y\}$. Moreover, $\ell_N(g') = m =  \ell_N(g'_1\cdots g'_{m-1}g''')$.
\qed

\bigskip\noindent
Now, combining Corollaries 3.7 and 3.8 and Lemmas 3.9 and 3.10, we get the following.

\bigskip\noindent
{\bf Proposition 3.11.}
{\it Let $M_1$, $M_2$ be preGarside monoids, and let $N$ be a common parabolic submonoid. Then the amalgamated product $M = M_1\!*_N\!M_2$ is a preGarside monoid.}
\qed

\bigskip\noindent
So, the amalgamated product of two preGarside monoids above a common parabolic submonoid is again preGarside. Moreover, the parabolic submonoids of the amalgamated product are as follows.

\bigskip\noindent
{\bf Proposition 3.12.}
{\it Let $M_1$, $M_2$ be preGarside monoids, and let $N$ be a common parabolic submonoid. Set $M = M_1\!*_N\!M_2$.
\begin{itemize}
\item[(1)]
If $M'_1$ and $M'_2$ are parabolic submonoids of $M_1$ and $M_2$, respectively, such that $M'_1\cap N = M'_2\cap N$, then $M'_1\!*_{N'}\!M'_2$ is (canonically isomorphic to) a parabolic submonoid of $M$, where $N' = M'_1\cap N$. In particular, $M_1$, $M_2$, and $N$ are parabolic submonoids of $M$.
\item[(2)]
If $M'$ is a parabolic submonoid of $M$, then there exist parabolic submonoids $M_1'$, $M_2'$, and $N'$, of $M_1$, $M_2$, and $N$, respectively, such that $M_1' \cap N = M_2' \cap N = N'$, and $M'$ is equal (isomorphic) to $M_1'\!*_{N'}\!M_2'$. 
\end{itemize}}

\bigskip\noindent
{\bf Proof.}
{\it Proof of (1).}
Let $M'_1$ and $M'_2$ be parabolic submonoids of $M_1$ and $M_2$, respectively, such that $M'_1\cap N = M'_2\cap N$. Set  $M' = M'_1\!*_{N'}\!M'_2$,  where $N' = M'_1\cap N$. If $g\in M'_i$ is such that $gN'$ is a minimal $\RR_{N'}$-class in $M'_i$, then $gN$ is a minimal $\RR_{N}$-class in $M_i$. Indeed, if $g = g_1h$ with $h$ in $N$, then $g_1 \in M_i'$ and $h \in N'$ because $M_i'$ is special. Hence, The canonical morphism from $M'$ to $M$ sends amalgam decompositions to amalgam decompositions, thus it is injective. 

\bigskip\noindent
Now, we prove that $M'$ is special. Let $g \in M'$ and $g',g'' \in M$ be such that $g=g'g''$. Set $[g]_N = (g_1, \dots, g_p, h)$, $[g']_N = (g_1', \dots, g_m', h')$, and $[g'']_N = (g_1'', \dots, g_n'',h'')$. Note that, since the embedding of $M'$ into $M$ sends amalgam decompositions to amalgam decompositions, we have $g_1, \dots, g_p \in M_1' \cup M_2'$ and $h \in N'$. Assume $n \ge 1$, $m \ge 1$, and $\varepsilon_m([g']_N) = \varepsilon_1([g'']_N)$. By Lemma~3.6, it follows that $g_1' = g_1, \dots, g_{m-1}' = g_{m-1}$, and
\[
[g_m'h'g'']_N=(g_m, g_{m+1}, \dots, g_p,h)\,.
\]
In particular, $g_1', \dots, g_{m-1}' \in M_1' \cup M_2'$. By Lemma 3.5, there exist $\tilde h_1, \dots, \tilde h_{n-1} \in N$ such that
\begin{gather*}
[g_m'h'g_1'']_N = (g_m, \tilde h_1)\,,\ [\tilde h_{i-1} g_i'']_N=(g_{m-1+i}, \tilde h_i) \text{ for } 2 \le i \le n-1\,,\\ 
\text{and } [\tilde h_{n-1} g_n'' h'']_N=( g_{m+n-1},h)\,.
\end{gather*}
Since $g_{m+n-1} h \in M_1' \cup M_2'$ and the submonoids $M_1'$ and $M_2'$ are special in $M_1$ and $M_2$, respectively, the equality $\tilde h_{n-1} g_n'' h'' = g_{m+n-1} h$ implies that $\tilde h_{n-1}, h'' \in N'$ and $g_n'' \in M_1' \cup M_2'$. Then, using the equality $\tilde h_{i-1} g_i'' = g_{m+i-1} \tilde h_i$, we prove by induction on $n-i$ that $\tilde h_{i-1} \in N'$ and $g_i'' \in M_1' \cup M_2'$ for $2 \le i \le n-1$. Finally, the equality $g_m' h' g_1'' = g_m \tilde h_1$ implies that $g_m', g_1'' \in M_1' \cup M_2'$ and $h' \in N'$. So, $g',g'' \in M'$. It is easily proved in the same way that $g',g'' \in M'$ if either $n=0$, or $m=0$, or $\varepsilon_m([g']_N) \neq \varepsilon_1([g'']_N)$.

\bigskip\noindent
Now, we take $g,g' \in M'$ such that $g \vee_{\!\scriptscriptstyle L} g'$ exists, and turn to prove that $g \vee_{\!\scriptscriptstyle L} g' \in M'$. Set $k = \ell_N(g)$ and $m=\ell_N(g')$, and suppose $k \le m$. We argue by induction on $m$, following the construction of $g \vee_{\!\scriptscriptstyle L} g'$ made in the proof of Lemma 3.10. If $m\le 1$, then there exists $i \in \{1,2\}$ such that $g,g' \in M_i'$. Then, since $M_i'$ is parabolic in $M_i$, we have $g \vee_{\!\scriptscriptstyle L} g' \in M_i'$. Now, assume $m \ge 2$ plus the inductive hypothesis. Set $[g]_N = (g_1, \dots, g_k,h)$ and $[g']_N = (g_1', \dots, g_m',h')$. Note that $g_1, \dots, g_k, g_1', \dots, g_m' \in M_1' \cup M_2'$ and $h,h' \in N'$. If $k \ge 2$, then $g_1=g_1', \dots, g_{k-1} = g_{k-1}'$, and
\[
g \vee_{\!\scriptscriptstyle L} g' = g_1 \cdots g_{k-1} ((g_kh) \vee_{\!\scriptscriptstyle L} (g_k' \cdots g_m'h'))\,.
\]
By the induction hypothesis we have $((g_kh) \vee_{\!\scriptscriptstyle L} (g_k' \cdots g_m'h')) \in M'$, thus $g \vee_{\!\scriptscriptstyle L} g' \in M'$. Hence, we may assume that $k \le 1$. Set $g'' = g \vee_{\!\scriptscriptstyle L} (g_1' \cdots g_{m-1}')$. By the induction hypothesis we have $g'' \in M'$. Moreover, as pointed out in the proof of Lemma 3.10, there exists $\tilde g \in M_1 \cup M_2$ such that $g'' = g_1' \cdots g_{m-1}' \tilde g$. Note that $\tilde g \in M'$ since $M'$ is special. Set $g''' = \tilde g \vee_{\!\scriptscriptstyle L} (g_m'h')$. Then, by the case $m=1$ treated above, $g''' \in M'$. On the other hand, it is shown in the proof of Lemma~3.10 that $g \vee_{\!\scriptscriptstyle L} g' = g_1' \cdots g_{m-1}' g'''$. So, $g \vee_{\!\scriptscriptstyle L} g' \in M'$.

\bigskip\noindent
{\it Proof of (2).}
Let $M'$ be a parabolic submonoid of $M$. Clearly the monoids $M'_1 = M_1 \cap M'$, $M'_2 = M_2 \cap M'$ and $N'=N \cap M'$ are parabolic in $M_1$, $M_2$ and $N$, respectively, and $N'=M_1' \cap N=M_2' \cap N$. Let $M'' = M'_1\!*_{N'}\!M'_2$. We claim that $M''$ is isomorphic to $M'$. Indeed, the image of $M''$ in $M$ by the canonical morphism is clearly included in $M'$. Conversely, if $g$ lies in $M'$ and $(g_1,\dots, g_n,h)$ is its amalgam decomposition in $M$, then each term $g_1,\dots, g_n$ and $h$ belongs to $M'$, because $M'$ is special, thus lies in $M_1' \cup M_2'$. So, $g \in M''$.
\qed


\section{PreGarside groups of FC type}

Now, thanks to the results of the previous section, mainly Proposition 3.11, we are ready to introduce the main definition of the paper. 

\bigskip\noindent
{\bf Definition.}
The family of \emph{preGarside monoids of FC type} is the smallest family of preGarside monoids that contains all Garside monoids and which is closed under amalgamation above parabolic submonoids. A preGarside group $G(M)$ is \emph{of FC type} if $M$ is.

\bigskip\noindent
As pointed out in the introduction, our goal in this section is to study preGarside groups of FC type. But, we need first to understand minimal coset representatives in Garside groups. This is the objective of the following subsection. 

\subsection{Minimal coset representatives in Garside groups}

Altobelli proved in~\cite{Alt} that, for each parabolic subgroup~$H$ of a spherical type Artin-Tits group~$G$, each left class $gH$ has a distinguished representative element~$m_H(g)$ that is minimal among the elements of~$gH$ for some partial order~$\leq_H$. Here we extend Altobelli's results to the wider context of Garside groups, with some new arguments and simplifications.

\bigskip\noindent
Throughout the subsection we assume $M$ is a Garside monoid with a Garside element $\Delta$, and $N$ is a parabolic submonoid of $M$. Recall from Lemma 2.4 that there is a Garside elements $\Delta_N$ of $N$ such that $\Div (\Delta) \cap N = \Div( \Delta_N)$. We start with two technical lemmas.

\bigskip\noindent
{\bf Lemma 4.1.}
{\it  Let $h_1,h_2$ belong to $N$ and $g$ lie in $M$. If $h_2 \wedge_{\!\scriptscriptstyle L} h_1 \neq h_2 \wedge_{\!\scriptscriptstyle L} (h_1g)$, then $\Delta_N \wedge_{\!\scriptscriptstyle L} g\neq 1$.}

\bigskip\noindent
{\bf Proof.}
Let $h_3$ be in $N$ such that $h_1 \vee_{\!\scriptscriptstyle L} (h_2 \wedge_{\!\scriptscriptstyle L} (h_1g)) = h_1h_3$. Both, $(h_2 \wedge_{\!\scriptscriptstyle L} (h_1g))$ and $h_1$, lie in $N$ and left divide $h_1g$, therefore $h_3$ belongs to $N$ and left divides $g$. If $h_3 = 1$, then $h_2 \wedge_{\!\scriptscriptstyle L} (h_1g)$ left divides both, $h_1$ and $h_2$, therefore it is equal to $h_2 \wedge_{\!\scriptscriptstyle L} h_1$. If $h_3 \neq 1$, then $h_3 \wedge_{\!\scriptscriptstyle L} \Delta_N$ is a common left divisor of $g$ and $\Delta_N$ different from $1$, hence $\Delta_N \wedge_{\!\scriptscriptstyle L} g\neq 1$.
\qed

\bigskip\noindent
Recall that a \emph{$\Delta$-simple} element is a factor of the Garside element $\Delta$, that is, an element of $\Div(\Delta)$. Throughout the subsection, for $g$ in $G(M)$, we denote by $|g|$ the smallest non-negative integer $k$ such that $g$ can be decomposed as a product of $k$ $\Delta$-simple elements and their inverses. Since here $\Delta$ is fixed, this does not induce confusion. We recall that a normal form is geodesic. In other words, if $g = ab^{-1}$ is in normal form, then $|g| = |a|+|b|$ (see \cite{DePa}).

\bigskip\noindent
{\bf Lemma 4.2.}
{\it Let $a$ and $b$ belong to $M$.
\begin{itemize}
\item[(1)]
The increasing sequence $(b \wedge_{\!\scriptscriptstyle L} \Delta_N^n)_{n\geq 0}$ stabilizes for $n\geq |b|$.
\item[(2)]
The increasing sequence $(a \wedge_{\!\scriptscriptstyle R} \Delta_N^nb)_{n\geq 0}$ stabilizes for $n\geq |a|$.
\end{itemize}}

\bigskip\noindent
{\bf Proof.}
By symmetry between left and right divisibilities, it suffices to prove (2). The sequence $(a\wedge_{\!\scriptscriptstyle R} \Delta_N^n b)_{n\geq 0}$ is bounded by $a$ for right divisibility, therefore it has to stabilize. Let $m$ be minimal such that $(a \wedge_{\!\scriptscriptstyle R} \Delta_N^nb)_{n\geq 0}$ stabilizes for $n\geq m$. We assume $m\geq 1$ and $a\neq 1$, otherwise there is nothing to prove. For short, we set $k = |a\wedge_{\!\scriptscriptstyle R} \Delta_N^mb|$. We are going to prove by induction on $k$ that $m\leq k$. This will prove (2), as $k \leq |a|$, since $a\wedge_{\!\scriptscriptstyle R} \Delta_N^mb$ right divides $a$. If $k = 0$, that is, $a\wedge_{\!\scriptscriptstyle R} \Delta_N^mb = 1$, the result follows from the fact that $(a\wedge_{\!\scriptscriptstyle R} \Delta_N^nb)_{n\geq 0}$ is increasing. Assume $k\geq 1$ plus the induction hypothesis. Denote by $c_k\cdots c_1$ the right greedy normal form of $a\wedge_{\!\scriptscriptstyle R} \Delta_N^mb$. Then there exists $a_1$ in $M$ such that $a = a_1c_1$. Moreover, 
\[
c_1 = \Delta \wedge_{\!\scriptscriptstyle R} (a \wedge_{\!\scriptscriptstyle R} \Delta_N^mb) = a\wedge_{\!\scriptscriptstyle R}(\Delta \wedge_{\!\scriptscriptstyle R} \Delta_N^mb) = a\wedge_{\!\scriptscriptstyle R} (\Delta \wedge_{\!\scriptscriptstyle R} \Delta_N b)\,.
\]
The last equality follows from the fact that $\Delta \wedge_{\!\scriptscriptstyle R} cd = \Delta \wedge_{\!\scriptscriptstyle R}((\Delta \wedge_{\!\scriptscriptstyle R} c)d)$ for all $c,d \in M$ (see \cite{DePa}). Therefore, there exists $b_1$ in $M$ such that $\Delta_N b = b_1c_1$, and $a\wedge_{\!\scriptscriptstyle R} \Delta_N^{n+1}b = (a_1\wedge_{\!\scriptscriptstyle R} \Delta_N^n b_1)c_1$ for every non-negative integer $n$. In particular, $(a_1 \wedge_{\!\scriptscriptstyle R} \Delta_N^n b_1)_{n\geq 0}$ stabilizes for $n \ge m-1$, and the right greedy normal form of $a_1\wedge_{\!\scriptscriptstyle R} \Delta_N^{m-1}b_1$ is $c_k \cdots c_2$. Applying the induction hypothesis, we get $m-1\leq k-1$ and we are done.
\qed

\bigskip\noindent
{\bf Definition.}
\begin{itemize}
\item[(1)]
For $g$ in $M$, we set
\[
M_N(g) = (g\wedge_{\!\scriptscriptstyle L} \Delta_N^{|g|})^{-1}g\,.
\]
\item[(2)]
We define the binary relation $\leq_N$ on $G(M)$ in the following way. Let $g_1 = a_1b_1^{-1}$ and $g_2 = a_2b_2^{-1}$ belong to $G(M)$ and be in normal form. We declare that $g_1\leq_N g_2$ if there exist $h_1,h_2$ in $N$ and $a$ in $M$ such that
\[
a_2 = a_1a\,,\ h_2b_2 = h_1b_1a\,,\  h_2\wedge_{\!\scriptscriptstyle L} h_1 = h_2 \wedge_{\!\scriptscriptstyle L} (h_1b_1)\,.
\]
\item[(3)]
We define the binary relation $\leq$ on $G(N)$ setting $h_1\leq h_2$ if $h_2h_1^{-1}$ belongs to $N$.
\item[(4)]
For $g$ in $G(M)$ with normal form $ab^{-1}$, we define the map $\varphi_g : N\to G(M)$ setting 
\[
\varphi_g(h) = a\, (a\wedge_{\!\scriptscriptstyle R} (hb))^{-1} \left(M_N\left(hb(a\wedge_{\!\scriptscriptstyle R} (hb))^{-1}\right)\right)^{-1}\,.
\]
\end{itemize}

\bigskip\noindent
{\bf Remark.}
\begin{itemize}
\item[(1)]
Let $ab^{-1}$ be the normal form of $g$. Set $c=a (a\wedge_{\!\scriptscriptstyle R}(hb))^{-1}$ and $d =  M_N(hb(a\wedge_{\!\scriptscriptstyle R} (hb))^{-1})$. It is easily checked that $c,d$ belong to $M$ and $c\wedge_{\!\scriptscriptstyle R} d = 1$. So, the formula in (4) provides the normal form for $\varphi_g(h)$.
\item[(2)]
The defining formula of $\varphi_g$ is quite ugly, but it is very easy to explain what this map does: starting with $h$, put $gh^{-1}$ in normal form $ef^{-1}$; then remove from $f^{-1}$ the inverse of the greatest left divisor of $f$ that lies in $N$. What remains is $\varphi_g(h)$. 
\end{itemize}

\bigskip\noindent
In order to prove Theorem 4.4, we need the following.

\bigskip\noindent
{\bf Lemma 4.3.}
{\it For $g$ in $G(M)$ and $h$ in $N$, one has $\varphi_g(h) = \varphi_{gh^{-1}}(1)$.}

\bigskip\noindent
{\bf Proof.}
It is easy to see that, for $a,b$ in $M$, even when $ab^{-1}$ is not a normal form, we have 
\[
\varphi_{ab^{-1}}(h) = a\, (a\wedge_{\!\scriptscriptstyle R} (hb))^{-1} \left(M_N\left(hb(a\wedge_{\!\scriptscriptstyle R} (hb)\right)^{-1})\right)^{-1}\,.
\]
In particular, for $g = ab^{-1}$ we get $\varphi_{gh^{-1}}(1) = \varphi_{a(hb)^{-1}}(1) = \varphi_g(h)$.
\qed

\bigskip\noindent
{\bf Theorem 4.4.}
{\it \begin{itemize}
\item[(1)]
The relation $\leq_N$ is a partial order on $G(M)$, and, if $g_1,g_2\in G(M)$ are such that $g_1\leq_N g_2$, then $g_1G(N) = g_2G(N)$.
\item[(2)]
The relation $\leq$ is a partial order on $G(N)$.
\item[(3)]
For every $g$, the map $\varphi_g$ is decreasing from $(N,\leq)$ to $(G(M),\leq_N)$, that is, $\varphi_g(h_2)\leq_N\varphi_g(h_1)$ if $h_1,h_2\in N$ are such that $h_1\leq h_2$. Moreover, for every $h$ in $N$, $\varphi_g(h)\leq_N gh^{-1}$.
\item[(4)]
For every $g$, the left coset $gG(N)$ has a unique minimal element $m_N(g)$ for the partial order $\leq_N$. Furthermore, $m_N(g) = \varphi_g(\Delta_N^{|g|})$.   
\end{itemize}}

\bigskip\noindent
{\bf Proof.}
{\it Proof of (1).}
The relation $\leq_N$ is clearly reflexive. Assume $g_1\leq_N g_2\leq_N g_1$, where $g_1 = a_1b_1^{-1}$ and $g_2 = a_2b_2^{-1}$ are in normal form. There exist $h_1,h_2$ in $N$ and $a$ in $M$ such that $a_2 = a_1a$, $h_2b_2 = h_1b_1a$, and $h_2\wedge_{\!\scriptscriptstyle L} h_1 = h_2 \wedge_{\!\scriptscriptstyle L} (h_1b_1)$. There exist also $h'_1,h'_2$ in $N$ and $a'$ in $M$ such that $a_1 = a_2a'$, $h'_2b_1 = h'_1b_2a'$, and $h'_2\wedge_{\!\scriptscriptstyle L} h'_1 = h'_2 \wedge_{\!\scriptscriptstyle L} (h'_1b_2)$. We can assume without restriction that $h_2\wedge_{\!\scriptscriptstyle L} h_1 = h'_2\wedge_{\!\scriptscriptstyle L} h'_1  =1$.  We get $a_2 = a_2a'a$. By cancellativity and atomicity it follows that $a = a' = 1$ and $a_1 = a_2$. Therefore $h_2b_2 = h_1b_1$ and $h_2 \wedge_{\!\scriptscriptstyle L} (h_1b_1) = 1$, which imposes $h_2 = 1$. Similarly, $h'_2b_1 = h'_1b_2$, $h'_2 \wedge_{\!\scriptscriptstyle L} (h'_1b_2) = 1$, and $h'_2 = 1$. We get $b_2 =h_1h'_1b_2$ and, by cancellativity and atomicity, $h'_1 = h_1 = 1$. Hence $b_1 = b_2$. So, the relation $\leq_N$ is anti-symmetric.

\bigskip\noindent
Now, assume $g_1\leq_N g_2$ and $g_2\leq_N g_3$, where $g_1 = a_1b_1^{-1}$, $g_2 = a_2b_2^{-1}$, and $g_3 = a_3b_3^{-1}$ are in normal form. There exist $h_1,h_2$ in $N$ and $a$ in $M$ such that $a_2 = a_1a$, $h_2b_2 = h_1b_1a$ and $h_2\wedge_{\!\scriptscriptstyle L} h_1 = h_2 \wedge_{\!\scriptscriptstyle L} (h_1b_1)$. There exist also $h'_2,h_3$ in $N$ and $a'$ in $M$ such that $a_3 = a_2a'$, $h_3b_3 = h'_2b_2a'$, and $h_3\wedge_{\!\scriptscriptstyle L} h'_2 = h_3 \wedge_{\!\scriptscriptstyle L} (h'_2b_2)$. As above, we assume without restriction that $h_2\wedge_{\!\scriptscriptstyle L} h_1 = h_3\wedge_{\!\scriptscriptstyle L} h'_2  =1$. We have $a_3 = a_1aa'$ and 
\[
((h_2\vee_{\!\scriptscriptstyle R} h'_2){h'_2}^{-1})h_3b_3 = (h_2\vee_{\!\scriptscriptstyle R} h'_2)b_2a' = ((h_2\vee_{\!\scriptscriptstyle R} h'_2){h_2}^{-1})h_1b_1aa'\,.
\]
Since $N$ is a parabolic submonoid, the elements $h'_3 = (h_2\vee_{\!\scriptscriptstyle R} h'_2){h'_2}^{-1}h_3$ and $h'_1 = (h_2\vee_{\!\scriptscriptstyle R} h'_2){h_2}^{-1}h_1$ belong to $N$. With this notation, we have $h'_3b_3 = h'_1b_1aa'$. It remains to show that $h'_3\wedge_{\!\scriptscriptstyle L}h'_1 = h'_3 \wedge_{\!\scriptscriptstyle L} (h'_1b_1)$. Let $c=h_3' \wedge_{\!\scriptscriptstyle L} (h_1'b_1)$. Write $h'_3 = cx$ and $h'_1b_1 = cy$. By cancellativity, we get $xb_3 = yaa'$. Since we assume $h_3\wedge_{\!\scriptscriptstyle L} h'_2b_2  =1$, we have $h_3b_3 = h'_2b_2a' = b_3\vee_{\!\scriptscriptstyle R} a'$. Therefore, there exists $z$ in $M$ such that $x = zh_3$ and $ya = zh'_2b_2$. As before, $h_2b_2 = h_1b_1a = b_2\vee_{\!\scriptscriptstyle R} a$, because $h_2\wedge_{\!\scriptscriptstyle L} h_1b_1  =1$, and there exists $t$ in $M$ such that $y = th_1b_1$. Hence $h'_1b_1 = cy = cth_1b_1$, and $c$ left divides $h'_1$ by cancellativitiy. So, $h'_3\wedge_{\!\scriptscriptstyle L} h'_1 = h'_3 \wedge_{\!\scriptscriptstyle L} (h'_1b_1)$. 

\bigskip\noindent
Finally, if $g_1\leq_N g_2$, then $g_2 = g_1(h_1^{-1}h_2)$, with the above used notations. Therefore, $g_1G(N) = g_2G(N)$.

\bigskip\noindent
{\it Proof of (2).}
Left to the reader.

\bigskip\noindent
{\it Proof of (3).}
Assume $h_2$ belongs to $N$. We have $\varphi_g(1)= a\, M_N(b)^{-1}$ and, by definition of $M_N(b)$, there exists $h$ in $N$ such that $b = h\,M_N(b)$. Similarly, there exists $h'$ in $N$ such that $h'M_N(h_2b(a\wedge_{\!\scriptscriptstyle R} (h_2b))^{-1}) = h_2b(a\wedge_{\!\scriptscriptstyle R} (h_2b))^{-1}$. Therefore,  
\[
h'M_N(h_2b(a \wedge_{\!\scriptscriptstyle R} (h_2b))^{-1})(a\wedge_{\!\scriptscriptstyle R} (h_2b)) = h_2hM_N(b)\,.
\]
Moreover, by Lemma 4.2, for $b$ in $M$, we have $M_N(b)\wedge_{\!\scriptscriptstyle L} \Delta_N = 1$.  Applying Lemma 4.1 we get 
\[
(h_2h)\wedge_{\!\scriptscriptstyle L} h'M_N(h_2b(a\wedge_{\!\scriptscriptstyle R} (h_2b))^{-1}) = (h_2h)\wedge_{\!\scriptscriptstyle L} h'\,.
\]
Hence, $\varphi_g(h_2)\leq_N\varphi_g(1)$. Moreover, if we assume $h_1\leq h_2$ in $N$ and write $h_2 = h_3h_1$, we have 
\[
\varphi_g(h_2) = \varphi_g(h_3h_1) = \varphi_{gh^{-1}_1}(h_3)\leq_N \varphi_{gh_1^{-1}}(1) = \varphi_g(h_1)\,.
\]

\bigskip\noindent
Now, the normal form  of $\varphi_g(1)$ is $a\,M_N(b)^{-1}$. But, by definition, there exists $h$ in $N$ such that $b = h\,M_N(b)$. This implies that $\varphi_g(1)\leq_N g$. Thus, for $h\in N$, we get $\varphi_g(h) = \varphi_{gh^{-1}}(1)\leq_N  gh^{-1}$.

\bigskip\noindent
{\it Proof of (4).} 
The coset $gG(N)$ contains minimal elements for $\leq_N$ by atomicity of $M$: if $(g_n)$ is a decreasing sequence for $\leq_N$ and $a_nb_n^{-1}$ is the normal form of $g_n$, then $a_{n+1}$ left-divides $a_n$ and the sequence $(a_n)$ has to stabilize. This implies that, for $n$ large enough, the sequence $(b_n)$ is decreasing for right divisibility. Therefore the sequence $b_n$ has to stabilize, too. Now, assume $gG(N) = g'G(N)$. There exists $h$ in $G(N)$ such that $gh = g'$. Let $h_1,h_2$ lie in $N$ such that $h = h^{-1}_1h_2$. Then $gh_1^{-1} = g'h_2^{-1}$. By (3), $\varphi_{gh^{-1}_1}(1) = \varphi_g(h_1)\leq_N \varphi_g(1)\leq_N g$ and  $\varphi_{g'h^{-1}_2}(1) = \varphi_{g'}(h_2)\leq_N \varphi_{g'}(1)\leq_N g'$. Assume $g'$ is minimal. Then $g' = \varphi_{gh^{-1}_1}(1) \leq_N g$. In particular, if $g$ is also minimal, then $g = g'$. Therefore, $gG(N)$ contains a unique minimal element, $m_N(g)$, for $\leq_N$, and there exists $h_1$ in $N$ such that $m_N(g) = \varphi_g(h_1)$. But, there exists a positive integer $k$ such that $h_1\leq \Delta_N^k$. Still by (3) and by minimality of $m_N(g)$, this implies that $\varphi_g(\Delta_N^k) = m_N(g)$. It remains to prove that the decreasing sequence $(\varphi_g(\Delta_N^k))_{k\geq 0}$ stabilizes at $k = |g|$. Assume that $g = ab^{-1}$ is in normal form and denote by $a_kb_k^{-1}$ the normal form of $\varphi_{g}(\Delta_N^k)$. It follows from the definition of the map $\varphi_g$ that the equality $a_k = a_{k+1}$ implies $\varphi_{g}(\Delta_N^k) = \varphi_{g}(\Delta_N^{k+1})$. Now,  by Lemma 4.2, for $k \geq |a|$, one has $a\wedge_{\!\scriptscriptstyle R} (\Delta_N^kb) = a\wedge_{\!\scriptscriptstyle R} (\Delta_N^{|a|}b)$, which implies $a_k = a_{|a|}$, and $\varphi_g(\Delta_N^k) = \varphi_g(\Delta_N^{|a|})$. Since $|g|\geq |a|$, we conclude that $m_N(g) = \varphi_g(\Delta_N^{|g|})$. 
\qed

\bigskip\noindent
{\bf Remark.}
\begin{itemize}
\item[(1)] 
It follows from the definition of the function $m_N$ that, for every $g$ in $M$, the element $m_N(g)$ is in $M$, $m_N(g^{-1}) = M_N(g)^{-1}$, and $m_N(g)N$ is the minimal $\RR_N$-class in $M$ which contains~ $g$.
\item[(2)]
For every $g$ in $G(M)$, $gG(N) = m_N(g)G(N)$, thus $m_N(m_N(g)) = m_N(g)$.
\end{itemize}

\bigskip\noindent
{\bf Proposition 4.5.}
{\it Let $g_1,g_2$ belong to $G(M)$, let $K$ be a parabolic submonoid of $M$, and assume $g_2$ belongs to $G(K)$. If $g_1\leq_Ng_2$, then $g_1$ belongs to $G(K)$. In particular, $m_N(g_2)$ lies in $G(K)$.}

\bigskip\noindent
{\bf Proof.}
We keep the notations of the definition of $\le_N$. We can assume without restriction that $h_1\wedge_{\!\scriptscriptstyle L} h_2 = 1$. The elements $a_2$ and $b_2$ have to lie in $K$. Therefore, $a_1$ and $a$ lie in $K$, too. But $h_2b_2 = (h_1b_1)a = b_2\vee_{\!\scriptscriptstyle R} a$. This implies that $h_1,h_2$ and $b_1$ lie in $K$. Thus $g_1 = a_1b_1^{-1}$ belongs to $G(K)$.
\qed

\bigskip\noindent
{\bf Proposition 4.6.}
{\it Assume $M$ is finitely generated, and denote by $w\mapsto\overline{w}$ the canonical map from $(\Div(\Delta)^{\pm})^*$ onto $G(M)$. There is an algorithm that associates a word $m^*_N(w)$ in $(\Div(\Delta)^{\pm})^*$ to every word $w$ in $(\Div(\Delta)^{\pm})^*$ such that
\begin{itemize}
\item[(a)]
$m_N(\overline{w}) = \overline{m_N^*(w)}$;
\item[(b)]
if $\overline{w}G(N) =\overline{v}G(N)$, then $m_N^*(w) = m_N^*(v)$;
\item[(c)] 
if $w$ is an element of $(\Div(\Delta)^\pm)^\ast$, $K$ is a parabolic submonoid, and $\Delta_K$ is the Garside element of $K$ satisfying $\Div(\Delta) \cap K = \Div(\Delta_K)$, and if $G(K) \cap \overline{w}G(N) \neq \emptyset$, then $m_N^\ast(w)$ belongs to $(\Div(\Delta_K)^{\pm})^\ast$;
\item[(d)]
if $w$ belongs to $\Div(\Delta)^*$, then $m_N^*(w)$ belongs to $\Div(\Delta)^*$ and $\overline{m_N^*(w)}$ represents the minimal $\RR_N$-class which contains $\overline{w}$.
\end{itemize}}

\bigskip\noindent
{\bf Proof.}
As recalled in Proposition 2.2, every element $g$ in $G(M)$ has a unique normal form $ab^{-1}$, and the elements $a$ and $b$ have unique right greedy normal forms $(a_p,\dots,a_1)$ and $(b_q,\dots,b_1)$, respectively, where the terms belong to $\Div(\Delta)$. So, to each element $g$ in $G(M)$ is associated a unique expression $a_p\cdots a_1b_1^{-1}\cdots b_q^{-1}$ on $\Div(\Delta)^{\pm}$. Moreover, there is an algorithm that, given $w$ such that $\overline{w} = g$, computes the words $v_1 = a_p\cdots a_1$ and $v_2 = b_q\cdots b_1$ \cite{Deh,GoG}. For $w$ in $(\Div(\Delta)^{\pm})^*$, we denote by $m^*_N(w)$ the unique above expression $a_p\cdots a_1b_1^{-1}\cdots b_q^{-1}$ associated with $m_N(\overline{w})$.

\bigskip\noindent
There are algorithms that, given two words $w_1,w_2$ on $\Div(\Delta)$, compute representing words on $\Div(\Delta)$ of the elements $\overline{w_1} \vee_{\!\scriptscriptstyle R} \overline{w_2}$ and $\overline{w_1}\wedge_{\!\scriptscriptstyle R} \overline{w_2}$. Starting from $w$, one can compute two words $a,b$ on $\Div (\Delta)$ such that $\overline{w} = \overline{a} \overline{b}^{-1}$ and $\overline{a} \wedge_{\!\scriptscriptstyle R} \overline{b} = 1$. One can then compute a representing word $c$ of $\overline{a} \wedge_{\!\scriptscriptstyle R} (\Delta_N^{|\overline{w}|} \overline{b})$, and compute representing words $d$, $e$ of $\overline{a}\, \overline{c}^{-1}$, $\Delta_N^{|\overline{w}|} \overline{b} \overline{c}^{-1}$, respectively. Then one can compute a representing word $f$ of $\overline{e} \wedge_{\!\scriptscriptstyle L} \Delta_N^{| \overline{e}|}$ and a representing word $g$ of $\overline{f}^{-1} \overline{e}$. Finally, the word $dg^{-1}$ is a representing word of $m_N( \overline{w})$ and can therefore be used to compute $m_N^\ast (w)$. Now, (b) follows from Theorem 4.4\,(4), (c) follows from Proposition 4.5, and (d) follows from the remark given after Theorem 4.4.
\qed

\subsection{Algebraic properties of preGarside groups of FC type}

The aim of this subsection is to extend Proposition 2.5 to preGarside groups of FC type (see Theorem 4.10). We will also extend Proposition 4.6 in the sense that, given a preGarside group $G$ of FC type, and given a parabolic subgroup $H$ of $G$, every coset $gH$ admits some ``special'' representative (see Theorem 4.11). We start with some technical results that will be useful in the remainder. The following is classical in the subject (see \cite{Ser} for instance). 

\bigskip\noindent
{\bf Proposition 4.7.}
{\it Let $G = G_1\!*_H\! G_2$ be the amalgamated product of two groups $G_1$ and $G_2$ over $H$. Let $C_1$ and $C_2$ be transversals of $G_1/H$ and $G_2/H$, respectively, that contain $1$. For all $g$ in $G$ there exists a unique sequence $(g_1,\dots,g_n,h)$ such that
\begin{itemize} 
\item[(a)]
$g = g_1\cdots g_n h$;
\item[(b)] 
$h\in H$ and $g_i\in (C_1\cup C_2)\setminus \{1\}$ for all $i\in\{ 1,\cdots,n\}$;
\item[(c)]
$g_i\in C_1\iff g_{i+1}\in C_2$ for all $i\in \{1,\cdots, n-1\}$.
\qed
\end{itemize}}

\bigskip\noindent
As in the case of the amalgamated product of monoids considered in Section 3, the above sequence $(g_1,\dots,g_n,h)$ is called the \emph{amalgam normal form} of $g$ relative to the amalgamated product of groups.

\bigskip\noindent
{\bf Lemma 4.8.}
{\it Let $M_1$, $M_2$ be two preGarside monoids, and assume $N$ is a parabolic submonoid of both, $M_1$ and $M_2$. Set $M = M_1\!*_N\!M_2$. If the morphisms $G(N)\to G(M_1)$ and $G(N)\to G(M_2)$ are into, then the group $G(M)$ is equal (isomorphic) to $G(M_1)\!*_{G(N)}\!G(M_2)$.}

\bigskip\noindent
{\bf Proof.}
By Theorem 2.6, $G(M)$ and $G(M_1)\!*_{G(N)}\!G(M_2)$ have the same group presentation.
\qed

\bigskip\noindent
Now, recall Properties (P1), (P2), and (P3) given in the statement of Proposition 2.5.

\bigskip\noindent
{\bf Lemma 4.9.}
{\it Let $M_1$, $M_2$ be two preGarside monoids that satisfy Properties (P1), (P2), and (P3). Assume $N$ is a parabolic submonoid of both, $M_1$ and $M_2$. Set $M = M_1\!*_N\!M_2$. Assume $M'$ is a parabolic submonoid of $M$. Set $M'_1 = M'\cap M_1$, $M'_2 = M'\cap M_2$ and $N' = M'\cap N$. Then the subgroup of $G(M)$ generated by $M'$ is isomorphic to $G(M'_1)\!*_{G(N')}\!G(M'_2)$, that is, to $G(M')$.} 

\bigskip\noindent
{\bf Proof.}
By Proposition 3.12 and Lemma 4.8, we know that $G(M')$ is isomorphic to $G(M'_1)\!*_{G(N')}\!G(M'_2)$. Since $M_1$ and $M_2$ satisfy (P1) and (P2), $N$ also satisfies (P1) and (P2), and $G(M'_1)$, $G(M'_2)$, and $G(N')$ inject in $G(M_1)$, $G(M_2)$, and $G(N)$, respectively. Moreover, by (P3), one has $G(N)\cap G(M'_1) = G(N')$ in $G(M_1)$ and $G(N)\cap G(M'_2) = G(N')$ in $G(M_2)$. This implies that the canonical morphism from $G(M'_1)\!*_{G(N')}\!G(M'_2)$ to $G(M_1)\!*_{G(N)}\!G(M_2)$ is injective. Its image is clearly the subgroup of $G(M)$ generated by $M'$.
\qed

\bigskip\noindent
{\bf Definition.}
Let $M$ be a preGarside monoid. We say that a finite labelled binary rooted tree is a \emph{FC tree} for $M$ if 
\begin{itemize}
\item[(a)]
each node is labelled by a FC type preGarside monoid;
\item[(b)]
the root is labelled by $M$, and each leaf is labelled by a Garside monoid;
\item[(c)]
each monoid attached to an inner node is the amalgam product above a common parabolic monoid of the monoids attached to the two child of the node.   
\end{itemize}
Note that, by definition, there exists a FC tree $\T$ for $M$ if and only if $M$ is of FC type. Moreover, in this case, for each monoid $N$ associated with a node of $\T$, there is an injective morphism $\iota_{(\T,N)}: N\to M$.  

\bigskip\noindent
The following two theorems will be proved together. 

\bigskip\noindent
{\bf Theorem 4.10.}
{\it Let $M$ be a preGarside monoid of FC type.
\begin{itemize}
\item[(P1)]
The natural morphism $\iota : M \to G(M)$ is injective.
\item[(P2)]
Let $N$ be a parabolic submonoid of $M$. The parabolic subgroup of $G(M)$ generated by $N$ is isomorphic to $G(N)$, and we have $G(N) \cap M = N$.
\item[(P3)]
Let $N, N'$ be parabolic submonoids of $M$. Then $N \cap N'$ is a parabolic submonoid, and $G(N) \cap G(N') = G(N \cap N')$.
\item[(P4)]
$G(M)$ is torsion free.
\end{itemize}}

\bigskip\noindent
{\bf Theorem 4.11.}
{\it Let $M$ be a preGarside monoid of FC type, let $P$ be a parabolic submonoid of $M$, and let $\T$ be a FC tree for $M$. Then there exists a map $m_{\T,P}:G(M) \to G(M)$ satisfying the following properties.
\begin{itemize}
\item[(a)]
$m_{\T,P}(g)\, G(P) = g\, G(P)$ for all $g \in G(M)$, and, if $g'\, G(P) = g\,G(P)$, then $m_{\T,P}(g) = m_{\T,P}(g')$, for all $g,g' \in G(M)$.
\item[(b)]
Let $g$ be an element of $M$. Then $m_{\T,P}(g)$ lies in $M$ and represents the minimal $\RR_P$-class containing $g$.
\item[(c)]
Let $M'$ be a parabolic submonoid of $M$, and let $g \in G(M)$. If $g\,G(P) \cap G(M') \neq \emptyset$, then $m_{\T,P}(g) \in G(M')$.
\end{itemize}}

\bigskip\noindent
{\bf Proof.}
We choose a FC tree $\T$ for $M$, we denote by $n$ the number of leafs of $\T$, and we argue by induction on $n$.

\bigskip\noindent
Assume $n=1$, thus $M$ is a Garside monoid. Then $M$ satisfies Properties (P1)-(P4) of Theorem~4.10 by Proposition 2.5. Let $P$ be a parabolic submonoid of $M$. Set $m_{\T,P} = m_P$. Then $m_{\T,P}$ satisfies Property (a) of Theorem 4.11 by Theorem 4.4, it satisfies Property (b) by the definition itself of $m_P$ (see the remark before Proposition 4.5), and it satisfies Property (c) by Proposition 4.5.

\bigskip\noindent
Now, we assume $n \ge 2$, plus the induction hypothesis. Let $M_1, M_2$ be the children of $M$ relatively to the tree $\T$, and let $N=M_1 \cap M_2$. So, $M=M_1\! *_N \! M_2$. For $i=1,2$, we denote by $\T_i$ the full subtree of $\T$ rooted at $M_i$. Clearly, $\T_i$ is a FC tree for $M_i$, and the number of leafs of $\T_i$ is strictly less than $n$. Hence, we may apply the induction  hypothesis to $M_i$ and $\T_i$. On the other hand, the FC tree $\T_{i,N}$, obtained from $\T_i$ by replacing each monoid $M'$ attached to a node by $M' \cap N$, is a FC tree for $N$. So, we may also apply the induction hypothesis to $N$ and $\T_{i,N}$.

\bigskip\noindent
By the above, for $i=1,2$, there exists a map $m_{\T_i,N} : G(M_i) \to G(M_i)$ satisfying Properties (a), (b), (c) of Theorem 4.11. Set $C_i = \{m_{\T_i,N}(g) \mid g \in G(M_i)\}$, and denote by $T_i$ the set of representatives of minimal $\RR_N$-classes in $M_i$. Then, by (a), $C_i$ is a transversal of $G(M_i)/G(N)$ and, by (b), $T_i$ is a subset of $C_i$ (and therefore $1 \in C_i$). It follows that the natural map $M=M_1\! *_N\!M_2 \to G(M) = G(M_1)\! *_{G(N)} \! G(M_2)$ sends amalgam normal forms to amalgam normal forms, thus it is injective. In particular, this shows that $M$ satisfies Property (P1) of Theorem 4.10. This also shows the following.
\begin{itemize}
\item[($*$)]
If $g \in M$ and $(g_1, \dots, g_\ell, h)$ is the amalgam normal form of $g$ relative to the decomposition $G(M) = G(M_1)\! *_{G(N)}\! G(M_2)$, then $g_i \in M_1 \cup M_2$ for all $i \in \{1, \dots, \ell\}$, and $h \in N$.
\end{itemize}

\bigskip\noindent
Let $M'$ be a parabolic submonoid of $M$. Set $M_1' = M' \cap M_1$, $M_2' = M' \cap M_2$, and $N' = M' \cap N$. Then $M' = M_1' \! *_{N'} \! M_2'$ (see Proposition 3.12), and the parabolic subgroup of $G(M)$ generated by $M'$ is $G(M') = G(M_1') \! *_{G(N')} \! G(M_2')$ (see Proposition 4.9). Let $g \in G(M')$, and let $(g_1, \dots, g_\ell,h)$ be the amalgam normal form of $g$. For $i \in \{1, \dots, \ell\}$, we denote by $\varepsilon(i)$ the element of $\{1,2\}$ such that $g_i \in G(M_{\varepsilon(i)})$. Since $g \in G(M') = G(M_1') \! *_{G(N')} \! G(M_2')$, $g$ can be written in the form $g=g_1' \cdots g_\ell'$ with $g_i' \in G(M_{\varepsilon(i)}')$ for all $i \in \{1, \dots, \ell\}$. Moreover, there exist $h_0, h_1, \dots, h_\ell \in G(N)$ such that 
\[
g_ih_i=h_{i-1}g_i' \text{ for } i\in \{1, \dots, \ell\}\,,\ h_0=1\,, \text{ and } h_\ell=h\,.
\]
We have $g_1h_1 = g_1'  \in G(M_{\varepsilon(1)}')$. By (c) applied to $m_{\T_{\varepsilon(1)}, N}$, it follows that $g_1 \in G(M_{\varepsilon(1)}')$. Then, we also have
\[
h_1 \in G(N) \cap G(M_{\varepsilon(1)}') = G(N \cap M_{\varepsilon(1)}') = G(N') \quad(\text{by (P3)})\,.
\]
It is easily proved in the same way, with an induction on $i$, that $g_i \in G(M_{\varepsilon(i)}')$ and $h_i \in G(N')$ for all $i \in \{1, \dots, \ell\}$. In particular, $h=h_\ell \in G(N')$. Finally, we have the following. 
\begin{itemize}
\item[($**$)]
If $g \in G(M')$ and $(g_1, \dots, g_\ell, h)$ is the amalgam normal form of $g$ relative to the decomposition $G(M) = G(M_1)\! *_{G(N)}\! G(M_2)$, then $g_i \in G(M_1') \cup G(M_2')$ for all $i \in \{1, \dots, \ell\}$, and $h \in G(N')$.
\end{itemize}

\bigskip\noindent
Let $g \in G(M') \cap M$. Let $(g_1, \dots, g_\ell,h)$ be the amalgam normal form of $g$. For $i \in \{1, \dots, \ell\}$, we denote by $\varepsilon(i)$ the element of $\{1,2\}$ such that $g_i \in G(M_{\varepsilon(i)})$. By ($*$) we have $g_i \in M_{\varepsilon(i)}$ for all $i \in \{1, \dots, \ell\}$, and $h \in N$. By ($**$) we have $g_i \in G(M_{\varepsilon(i)}')$ for all $i \in \{1, \dots, \ell\}$, and $h \in G(N')$. By (P2) applied to $M_1$, $M_2$, and $N$, we have $G(M_{\varepsilon(i)}') \cap M_{\varepsilon(i)} = M_{\varepsilon(i)}'$ for all $i \in \{1, \dots, \ell\}$, and $G(N') \cap N = N'$. So, $g \in M'$. This shows that $M$ satisfies Property (P2) of Theorem 4.10.

\bigskip\noindent
Let $M', M''$ be parabolic submonoids of $M$. First, note that $M' \cap M''$ is parabolic by definition. Set $M_1' = M' \cap M_1$, $M_2' = M' \cap M_2$, $N' = M' \cap N$, $M_1'' = M'' \cap M_1$, $M_2'' = M'' \cap M_2$, and $N'' = M'' \cap N$. Let $g \in G(M') \cap G(M'')$. Let $(g_1, \dots, g_\ell,h)$ be the amalgam normal form of $g$. For $i \in \{1, \dots, \ell\}$, we denote by $\varepsilon(i)$ the element of $\{1,2\}$ such that $g_i \in G(M_{\varepsilon(i)})$. By ($**$) we have $g_i \in G(M_{\varepsilon(i)}') \cap G(M_{\varepsilon(i)}'')$ for all $i \in \{1, \dots, \ell\}$, and $h \in G(N') \cap G(N'')$. By (P3) applied to $M_1$, $M_2$, and $N$, we have $G(M_{\varepsilon(i)}') \cap G(M_{\varepsilon(i)}'') = G(M_{\varepsilon(i)}' \cap M_{\varepsilon(i)}'')$ for all $i \in \{1, \dots, \ell\}$, and $G(N') \cap G(N'') = G(N' \cap N'')$. This implies that $g \in G(M' \cap M'')$. The inclusion $G(M' \cap M'') \subset G(M') \cap G(M'')$ being trivial, this shows that $M$ satisfies Property~(P3) of Theorem 4.10.

\bigskip\noindent
Let $g$ be a finite order element of $G(M)$. Since $G=G(M_1)\! *_{G(N)} \! G(M_2)$, $g$ is conjugate to an element of either $G(M_1)$ or $G(M_2)$ (see \cite[p. 54]{Ser}). But, by the induction hypothesis, $G(M_1)$ and $G(M_2)$ are torsion free, thus $g=1$. This shows that $M$ satisfies (P4).

\bigskip\noindent
Now, we take a parabolic submonoid $P$ of $M$. We set $P_1 = M_1 \cap P$ and $P_2= M_2 \cap P$. Let $g \in G(M)$, and let $(g_1, \dots, g_\ell,h)$ be the amalgam normal form of $g$. For $i \in \{1, \dots, \ell\}$, we denote by $\varepsilon(i)$ the element of $\{1,2\}$ such that $g_i \in G(M_{\varepsilon(i)})$. We define $m_{\T,P}(g)$ by induction on $\ell$ as follows. Suppose $\ell=0$. Then $g \in G(N)$, and we set 
\[
m_{\T,P}(g) = m_{\T_{1,N}, P_1 \cap P_2} (g)\,.
\] 
Suppose $\ell=1$. Then $g \in G(M_{\varepsilon(1)})$, and we set
\[
m_{\T,P}(g) = m_{\T_{\varepsilon(1)},P_{\varepsilon(1)}}(g)\,.
\]
Suppose $\ell \ge 2$. If $m_{\T_{\varepsilon(\ell)},P_{\varepsilon(\ell)}} (g_\ell h) \not\in G(N)$, we set
\[
m_{\T,P} (g) = g_1 \cdots g_{\ell-1}\cdot m_{\T_{\varepsilon(\ell)},P_{\varepsilon(\ell)}}(g_\ell h)\,.
\]
If $m_{\T_{\varepsilon(\ell)},P_{\varepsilon(\ell)}} (g_\ell h) \in G(N)$, we set 
\[
g'=g_1 \cdots g_{\ell -1}\cdot m_{\T_{\varepsilon(\ell)},P_{\varepsilon(\ell)}}(g_\ell h) \quad \text{and} \quad m_{\T,P} (g) =  m_{\T,P} (g')\,.
\]

\bigskip\noindent
Let $g \in G(M)$, and let $(g_1, \dots, g_\ell,h)$ be the amalgam normal form of $g$. For $i \in \{1, \dots, \ell\}$, we denote by $\varepsilon(i)$ the element of $\{1,2\}$ such that $g_i \in G(M_{\varepsilon(i)})$. The number $\ell$ will be called \emph{amalgam length} of $g$, and it will be denoted by $|g|_a$. It is easily proved by induction on $|g|_a$ that
\[
m_{\T,P}(g)\,G(P) = g\,G(P)\,.
\]

\bigskip\noindent
We turn now to show that $|m_{\T,P}(g)|_a$ is minimal among the amalgam lengths of the elements of the coset $g\,G(P)$. Since $m_{\T,P}(g)\,G(P) = g\,G(P)$, we can assume that $g = m_{\T,P}(g)$. If $m_{\T,P}(g) \in G(N)$, then $| m_{\T,P}(g)|_a = 0$ is obviously minimal. We can therefore assume that $m_{\T,P}(g) \not \in G(N)$, thus, by construction, $g_\ell h = m_{\T_{\epsilon(\ell)}, P_{\epsilon(\ell)}} (g_\ell h) \not\in G(N)$. Let $u$ be in $G(P)$, and let $(u_1, \dots, u_k,v)$ be the amalgam normal form of $u$. If we had $|gu|_a < |g|_a$, then we would have $k\ge 1$, $u_1 \in G(P_{\varepsilon(\ell)})$ (by ($**$)) and $g_\ell hu_1 \in G(N)$, thus $(g_\ell h)G(P_{\varepsilon(\ell)}) \cap G(N) \neq \emptyset$, thus, by (c) of Theorem 4.11 applied to $m_{\T_{\varepsilon(\ell)}, P_{\varepsilon(\ell)}}$, we would have $g_\ell h = m_{\T_{\varepsilon(\ell)}, P_{\varepsilon(\ell)}}(g_\ell h) \in G(N)$: a contradiction.

\bigskip\noindent
Now, we take $g,g' \in G(M)$ such that $g\,G(P) = g'\, G(P)$, and we prove that $m_{\T,P} (g) = m_{\T,P} (g')$. By the above, we can assume that $g=m_{\T,P}(g)$, $g' = m_{\T,P }(g')$, and $|g|_a = |g'|_a$. Let $u \in G(P)$ such that $g'=gu$, and let $(g_1, \dots, g_\ell, h)$, $(g_1', \dots, g_\ell', h')$, and $(u_1, \dots, u_k,v)$ be the amalgam normal forms of $g$, $g'$, and $u$, respectively. Assume that $|g|_a = |g'|_a = 0$. Then we have $g=h$, $g'=h'$, and $u=v \in G(P_1 \cap P_2)$, thus $g G(P_1 \cap P_2) = g' G(P_1 \cap P_2)$. By the induction hypothesis, it follows that $g = m_{\T_{1,N},P_1 \cap P_2} (g) = m_{\T_{1,N},P_1 \cap P_2} (g')  = g'$. Assume that $|g|_a = |g'|_a > 0$. It is easily shown in the same way as before that $g_\ell h u_1 \not\in G(N)$, thus $k \le 1$ and $u\in G(P_{\varepsilon(\ell)})$ (by ($**$)). By the uniqueness of the amalgam normal form, it follows that $g_1 = g_1', \dots, g_{\ell -1} = g_{\ell -1}'$, and $g_\ell 'h' = g_\ell h u$. Since $u\in G(P_{\varepsilon(\ell)})$, we have $g_\ell h\, G(P_{\varepsilon(\ell)}) = g_\ell' h'\, G(P_{\varepsilon(\ell)})$, thus, by (a) of Theorem 4.11 applied to $m_{\T_{\varepsilon(\ell)}, P_{\varepsilon(\ell)}}$, 
\[
g_\ell h = m_{\T_{\varepsilon(\ell)}, P_{\varepsilon(\ell)}} (g_\ell h) = m_{\T_{\varepsilon(\ell)}, P_{\varepsilon(\ell)}} (g_\ell' h') = g_\ell' h'\,.
\]
So, $g=g'$.

\bigskip\noindent
It is easily shown by induction on $|g|_a$ that, if $g \in M$, then $m_{\T,P}(g)$ also belongs to $M$. 
We turn now to prove that, in that case, $m_{\T,P}(g)$ is the representative of the minimal $\RR_P$-class containing $g$. 
We assume without loss of generality that $g = m_{\T,P}(g)$.

\bigskip\noindent
Assume first that $|g|_a = 0$, that is, $g \in N$. 
Then $g = m_{\T_{1,N},P_1 \cap P_2}(g)$. 
Suppose $g$ is written $g=g'u$ with $g \in M$ and $u \in P$. 
Since $N$ is parabolic, we must have $g' \in N$ and $u \in P \cap N = P_1 \cap P_2$. 
Then, by the induction hypothesis (on the number of leaves of the FC tree), $u=1$ and $g'=g$.

\bigskip\noindent
Assume now that $|g|_a \ge 1$. 
Let $(g_1, \dots, g_\ell, h)$ be the amalgam normal form of $g$.
We have by construction $g_\ell h = m_{\T_{\varepsilon(\ell)},M_{\varepsilon(\ell)}}(g_\ell h)$.
Suppose $g$ is written $g=g'u$ with $g' \in M$ and $u \in P$.
By the above, $|g|_a$ is minimal among the amalgam lengths of the elements of the coset $g\,G(P)$.
By Lemma 3.6, it follows that the amalgam normal form of $g'$ is of the form $(g_1, \dots, g_{\ell-1}, g_\ell',h')$, that $u \in M_{\varepsilon(\ell)}$ (namely, $u \in P_{\varepsilon(\ell)}$), and that $g_\ell h = g_\ell'h'u$.
By the induction hypothesis, we conclude that $u=1$, thus $g'=g$.
So, $m_{\T,P}$ satisfies Property (b).

\bigskip\noindent
Let $M'$ be a parabolic submonoid of $M$. 
It is easily shown by induction on $|g|_a$ that, if $g \in G(M')$, then $m_{\T,P}(g) \in G(M')$. 
It follows that, if $g\,G(P) \cap G(M') \neq \emptyset$, then $m_{\T,P}(g) \in G(M')$, that is, $m_{\T,P}$ satisfies Property (c). 
Indeed, if $g' \in g\,G(P) \cap G(M')$, then, by the above, $m_{\T,P}(g) = m_{\T,P}(g') \in G(M')$.
\qed

\subsection{Combinatorial properties of FC type preGarside groups}

In this subsection we assume given a finite set $S$, two binary graphs $\Gamma_L$ and $\Gamma_R$ on $S$, a partial complement $f_L$ on $S$ based on $\Gamma_L$, and a partial complement $f_R$ on $S$ based on $\Gamma_R$, and we assume that $M=M_L(\Gamma_L,f_L) = M_R(\Gamma_R, f_R)$ is a FC type preGarside monoid. Recall that, thanks to Lemma 2.7, we can and we do assume that $S$ is the set of atoms of $M$. 

\bigskip\noindent
{\bf Remark.}
We cannot remove the assumption ``$M$ is a preGarside monoid of FC type'' because we do not know how to decide whether a monoid of the form $M=M_L(\Gamma_L,f_L) = M_R (\Gamma_R, f_R)$ is a preGarside monoid of FC type.

\bigskip\noindent
A direct consequence of Theorem 2.8 is the following.

\bigskip\noindent
{\bf Proposition 4.12.}
{\it Let $X_1, X_2$ be two non-empty subsets of $S$ such that $S=X_1 \cup X_2$, let $M_1, M_2$ be the submonoids of $M$ generated by $X_1, X_2$, respectively, and let $N=M_1 \cap M_2$. Then $M_1,M_2$ are parabolic submonoids of $M$ and $M=M_1\! *_N\! M_2$ if and only if the following two conditions hold.
\begin{itemize}
\item[(a)]
Let $i \in \{1,2\}$. For all $x,y \in X_i$, either $\{x,y\}$ is not an edge of $\Gamma_L$ or $f_L(x,y) \in X_i^*$, and either $\{x,y\}$ is not an edge of $\Gamma_R$ or $f_R(x,y) \in X_i^*$.
\item[(b)]
Let $i \in \{1,2 \}$. For all $x \in X_i$ and all $y \in S \setminus X_i$, the pair $\{x,y\}$ is an edge of neither $\Gamma_L$, nor $\Gamma_R$.
\qed
\end{itemize}}

\bigskip\noindent
The next result is easily proved from the above proposition. 

\bigskip\noindent
{\bf Corollary 4.13.}
{\it \begin{itemize}
\item[(1)]
$M$ is a Garside monoid if and only if $\Gamma_L = \Gamma_R= K_S$ is the complete graph on $S$.
\item[(2)]
There is an algorithm which determines a FC tree for $M$.
\item[(3)]
We have $\Gamma_L = \Gamma_R$.
\qed
\end{itemize}} 

\bigskip\noindent
{\bf Definition.}
If $M$ is a Garside monoid, we denote by $\SS(M)$ the set of divisors of the minimal Garside element $\delta=\delta_M$ of $M$. The elements of $\SS(M)$ are called \emph{simple elements} of $M$. Note that, if $N$ is a parabolic submonoid of $M$, then $\SS(N) \subset \SS(M)$. Let $M$ be a finitely generated preGarside monoid. Then $\SS(M)$ denotes the union of the simple elements of the spherical parabolic submonoids of $M$. Note that, if $N$ is a (spherical) parabolic submonoid of $M$, then $\AA(N) = \AA(M) \cap N$, thus there are finitely many spherical parabolic submonoids of $M$, therefore $\SS(M)$ is finite. As for Garside monoids, the elements of $\SS(M)$ are called \emph{simple elements}.

\bigskip\noindent
{\bf Remark.}
We do not necessarily have $\AA(M) \subset \SS(M)$ in general. For example, the monoid $\langle a,b,c \mid a^2=bc\rangle^+$ is a preGarside monoid whose unique spherical parabolic submonoid is $\{1\}$, hence $\SS(M) = \emptyset$ for this monoid while $\AA(M) = \{a,b,c\}$. However, it is easily seen that $\AA(M) \subset \SS(M)$ if $M$ is of FC type because it is so for Garside monoids.

\bigskip\noindent
Now, we come back to the hypothesis of the subsection: $M$ is a preGarside monoid of FC type given by two partial complements $f_L$ and $f_R$.

\bigskip\noindent
{\bf Lemma 4.14.}
{\it There exists an algorithm which determines $\SS(M)$, where each element $x \in \SS(M)$ is given by a word $a(x) \in S^*$.}

\bigskip\noindent
{\bf Proof.}
Let $X$ be a subset of $S$, and let $N$ be the submonoid of $M$ generated by $X$. Then, by Theorem 2.8 and Proposition 4.12, $N$ is a spherical parabolic submonoid if and only if the following hold.
\begin{itemize}
\item[(a)]
For all $x,y \in X$, $\{x,y\}$ is an edge of $\Gamma_L=\Gamma_R$, $f_L(x,y) \in X^*$, and $f_R(x,y) \in X^*$.
\item[(b)]
For all $x \in X$ and $y \in S \setminus X$, if $\{x,y\}$ is an edge of $\Gamma_L=\Gamma_R$, then $f_L(x,y) \not\in X^*$ and $f_R(y,x) \not\in X^*$.
\end{itemize}
In particular, there is an effective way of determining all spherical parabolic submonoids of $M$.

\bigskip\noindent
Now, suppose that $N$ is a spherical parabolic submonoid. Let $w \in X^*$, and let $g$ be the element of $N$ represented by $w$. Using any solution to the word problem in $N$ (see \cite{DePa}, \cite{Deh}, for example), we can determine all left and right factors of $g$, hence we can decide if $g$ is a Garside element. In order to calculate the minimal Garside element $\delta_N$ of $N$ as well as $\SS(N)=\Div (\delta_N)$, we apply this test to all words of length $1$, then to those of length $2$, an so on. We increase the length of the tested words until we obtain a Garside element, which should be the minimal one. 
\qed

\bigskip\noindent
{\bf Theorem 4.15.}
{\it Let $P$ be a parabolic submonoid of $M$, and let $\T$ be an FC tree for $M$. There exists a function $m^*_{\T,P}: \SS(M)^{\pm *} \to \SS(M)^{\pm *}$ satisfying the following properties.
\begin{itemize}
\item[(a)]
Let $w \in \SS(M)^{\pm *}$. Then $\overline{m^*_{\T,P} (w)} = m_{\T,P}( \overline{w})$.
\item[(b)]
Let $v,w \in \SS(M)^{\pm *}$. If $\overline{w}G(P) = \overline{v}G(P)$, then $m^*_{\T,P}(w) = m^*_{\T,P}(v)$.
\item[(c)]
Let $w \in \SS(M)^{\pm *}$. If $M'$ is a parabolic submonoid of $M$ and $G(M') \cap \overline{w} G(P) \neq \emptyset$, then $m^*_{\T,P}(w) \in (\SS(M)\cap M')^{\pm *}$.
\item[(d)]
Let $w \in \SS(M)^{\pm *}$. If $\overline{w} \in M$, then $m^*_{\T,P}(w) \in \SS(M)^*$.
\end{itemize}
Moreover, there is an algorithm which, given $w \in \SS(M)^{\pm *}$, determines $m^*_{\T,P}(w)$.}

\bigskip\noindent
{\bf Proof.}
We argue by induction on the number $n$ of leafs of $\T$. If $n=1$, then $M$ is a Garside monoid. In this case we set $m^*_{\T,P}=m^*_P$, and Properties (a), (b), (c) and (d) are satisfied by Proposition 4.6. So, we may assume that $n \ge 2$ plus the induction hypothesis.

\bigskip\noindent
Let $M_1, M_2$ be the children of $M$ relative to $\T$, and let $N=M_1 \cap M_2$. 
For $i=1,2$ we denote by $\T_i$ the full subtree of $\T$ rooted at $M_i$. 
On the other hand, we denote by $\T_{i,N}$ the FC tree for $N$ obtained from $\T_i$ by replacing each monoid $M'$ attached to a node by $M' \cap N$. 
We go back to the constructions and notations given in the proof of Theorems 4.10 and 4.11. 

\bigskip\noindent
Suppose $M'$ is a preGarside monoid of FC type, $P'$ is a parabolic submonoid, and $\T'$ is a FC tree for $M'$, and suppose that the number of leafs of $\T'$ is strictly less than $n$. Then, by the induction hypothesis, for all $g \in G(M')$, there exists a unique word $\omega_{\T',P'}(g)$ such that $\omega_{\T',P'}(g) = m^*_{\T',P'}(w)$ for all $w \in \SS(M')^{\pm *}$ such that $\overline{w} =g$. This word $\omega_{\T',P'}(g)$ will be used throughout the whole proof for $\T'=\T_1$ or $\T_2$ or $\T_{1,N}$ and $M'=M_1$ or $M_2$ or $N$, respectively.

\bigskip\noindent
Let $g \in G(M)$. Let $(g_1, \dots, g_\ell,h)$ be the amalgam normal form of $m_{\T,P}(g)$. For $i \in \{1, \dots, \ell\}$,  we denote by $\varepsilon(i)$ the element of $\{1,2\}$ such that $g_i \in M_{\varepsilon(i)}$. We set $u_i = \omega_{\T_{\varepsilon(i)},N} (g_i)$ for all $i \in \{1, \dots, \ell\}$, $v=\omega_{\T_{1,N},\{1\}}(h)$, and 
\[
\mu(g) = u_1 \cdots u_\ell v\,.
\]
For $w \in \SS(M)^{\pm *}$ we set
\[
m^*_{\T,P} (w) = \mu(\overline{w})\,.
\]

\bigskip\noindent
The fact that the function $m^*_{\T,P}$ satisfies Properties (a) and (b) follows from the construction of the function itself. Let $w \in \SS(M)^{\pm *}$, and let $M'$ be a parabolic submonoid of $M$ such that $G(M') \cap \overline{w} G(P) \neq \emptyset$. Set $M_1' = M_1 \cap M'$, $M_2' = M_2 \cap M'$, and $N' = N \cap M'$. Let $(g_1, \dots, g_\ell,h)$ be the amalgam normal form of $m_{\T,P}(\overline{w})$. By Theorem 4.11, we have $m_{\T,P}(\overline{w}) \in G(M')$ and, by Property~($**$) proved in the proof of Theorems 4.10 and 4.11, we have $g_i \in G(M_{\varepsilon(i)}')$ for all $i \in \{1, \dots, \ell\}$, and $h \in G(N')$. Then, by the induction hypothesis, $\omega_{\T_{\varepsilon(i)},P_{\varepsilon(i)}}(g_i) \in (\SS(M_{\varepsilon(i)}) \cap M_{\varepsilon(i)}')^{\pm *}$ for all $i \in \{1, \dots, \ell\}$, and $\omega_{\T_{1,N},\{1\}}(h) \in (\SS(N) \cap N')^{\pm *}$. This implies that $m^*_{\T,P}(w) \in (\SS(M) \cap M')^{\pm *}$.

\bigskip\noindent
Let $w \in \SS(M)^{\pm *}$ be such that $\overline{w} \in M$. Let $(g_1, \dots, g_\ell,h)$ be the amalgam normal form of $m_{\T,P}(\overline{w})$. By Theorem 4.11, we have $m_{\T,P}(\overline{w}) \in M$ and, by Property ($*$) proved in the proof of Theorems 4.10 and 4.11, we have $g_i \in M_{\varepsilon(i)}$ for all $i \in \{1, \dots, \ell\}$, and $h \in N$. Then, by the induction hypothesis, $\omega_{\T_{\varepsilon(i)},P_{\varepsilon(i)}}(g_i) \in \SS(M_{\varepsilon(i)})^*$ for all $i \in \{1, \dots, \ell\}$, and $\omega_{\T_{1,N},\{1\}}(h) \in \SS(N)^*$. This implies that $m^*_{\T,P}(w) \in \SS(M)^*$.

\bigskip\noindent
It remains to show that there is an algorithm which, given $w \in \SS(M)^{\pm *}$, determines $m^*_{\T,P}(w)$. Recall that, by hypothesis, the given generating set for $M$ is $S=\AA(M)$, and every element $x$ in $\SS(M)$ is given by a word $a(x) \in \AA(M)^*$. The map $a : \SS(M) \to \AA(M)^*$ induces a morphism $a^* : \SS(M)^{\pm *} \to \AA(M)^{\pm *}$ which will be useful in our construction. 

\bigskip\noindent
Define a {\it pre-expression} of length $\ell$ to be a pair of sequences
\[
W=((u_1, \dots, u_\ell,v),(\varepsilon(1), \dots, \varepsilon(\ell)))
\]
such that $\varepsilon(i) \in \{1,2\}$ and $u_i \in \SS(M_{\varepsilon(i)})^{\pm *}$ for all $i \in \{1, \dots, \ell\}$, and $v \in \SS(N)^{\pm *}$.

\bigskip\noindent
Let $W$ be a pre-expression. Suppose that $\varepsilon(i) = \varepsilon(i+1)$ for some $i \in \{1, \dots, \ell-1\}$. Set
\[
W' = ((u_1, \dots, u_iu_{i+1}, u_{i+2}, \dots,u_\ell),(\varepsilon(1), \dots, \varepsilon(i),\varepsilon(i+2), \dots, \varepsilon(\ell)))\,.
\]
Then $W'$ is called an elementary reduction of type I of $W$. 

\bigskip\noindent
Suppose $m^*_{\T_{\varepsilon(i)},N}(u_i)=1$, that is, $\overline{u_i} \in N$, for some $i \in \{1, \dots, \ell\}$. Set $u_i'= m^*_{\T_{\varepsilon(i)},\{1\}} (u_i) \in (\SS(M) \cap N)^{\pm *}$, and set $u_i'' = a^*(u_i') \in \AA(N)^{\pm *}$. Set 
\[
W'=\left\{ \begin{array}{ll}
((u_1, \dots, u_{i-1},u_i''u_{i+1}, \dots,u_\ell,v),(\varepsilon(1), \dots, \varepsilon(i-1), \varepsilon(i+1), \dots, \varepsilon(\ell))) & \text{if } i < \ell\\
((u_1, \dots, u_{\ell-1}, u_\ell''v),(\varepsilon(1), \dots, \varepsilon(\ell-1))) &\text{if } i=\ell
\end{array} \right.
\]
Then $W'$ is called an elementary reduction of type II of $W$. 

\bigskip\noindent
Suppose $m^*_{\T_{\varepsilon(i)},N}(u_i) \neq 1$ and $m^*_{\T_{\varepsilon(i)},N}(u_i) \neq u_i$, for some $i \in \{1, \dots, \ell\}$. Set  $u_i' = m^*_{\T_{\varepsilon(i)},N}(u_i)$ and $v_i = m^*_{\T_{\varepsilon(i)},N}(u_i)^{-1}\, u_i$. Note that $\overline{v_i} \in G(N)$. Set $v_i' = m^*_{\T_{\varepsilon(i)},\{1\}} (v_i) \in (\SS(M_{\varepsilon(i)}) \cap N)^{\pm *}$, and set $v_i'' = a^* (v_i') \in \AA(N)^{\pm *}$. Set
\[
W' = \left\{ \begin{array}{ll}
((u_1, \dots, u_i', v_i''u_{i+1}, \dots, u_\ell),(\varepsilon(1), \dots, \varepsilon(i), \varepsilon(i+1), \dots, \varepsilon(\ell))) &\text{if } i<\ell\\
((u_1, \dots, u_{\ell-1}, u_\ell',v_\ell''v),(\varepsilon(1), \dots, \varepsilon(\ell -1), \varepsilon(\ell))) &\text{if } i=\ell
\end{array} \right.
\]
Then $W'$ is called an elementary reduction of type III of $W$.

\bigskip\noindent
Suppose $m^*_{\T_{1,N},\{1\}}(v) \neq v$. Set $v' = m^*_{\T_{1,N},\{1\}}(v)$. Set 
\[
W' = ((u_1, \dots, u_\ell, v'),(\varepsilon(1), \dots, \varepsilon(\ell)))\,.
\]
Then $W'$ is called an elementary reduction of type IV of $W$.

\bigskip\noindent
Set $u_\ell'=m^*_{\T_{\varepsilon(\ell)},P_{\varepsilon(\ell)}}(u_\ell v)$, and suppose $m^*_{\T_{\varepsilon(\ell)},\{1\}} (u_\ell') \neq m^*_{\T_{\varepsilon(\ell)},\{1\}} (u_\ell v)$, that is, $\overline{u_\ell'} \neq \overline{u_\ell v}$. Set 
\[
W' = ((u_1, \dots, u_{\ell-1}, u_\ell',1),(\varepsilon(1), \dots, \varepsilon(\ell)))\,.
\]
Then $W'$ is called an elementary reduction of type V of $W$. 

\bigskip\noindent
A pre-expression $W'$ is called a {\it reduction} of $W$ if there is a finite sequence $W_0=W, W_1, \dots,W_p=W'$ of pre-expressions such that $W_i$ is an elementary reduction of $W_{i-1}$ for all $i \in \{1, \dots, p\}$. We say that a pre-expression $W$ is {\it reduced} if it has no elementary reductions. Observe that any sequence of elementary reductions is finite. Moreover, if $W=((u_1, \dots, u_\ell,v),(\varepsilon(1), \dots, \varepsilon(\ell)))$ is reduced, then $m^*_{\T,P}(u_1 \cdots u_\ell v) = u_1 \cdots u_\ell v$.

\bigskip\noindent
Now, let $w =s_1^{e_1} \cdots s_\ell ^{e_\ell} \in \AA(M)^{\pm *}$. For all $i \in \{1, \dots, \ell\}$ choose $\varepsilon(i) \in \{1,2\}$ such that $s_i \in \AA(M_{\varepsilon(i)})$. Let
\[
W_0 = ((s_1^{e_1}, \dots, s_\ell^{e_\ell},1),(\varepsilon(1), \dots, \varepsilon(\ell)))\,.
\]
Then a reduced reduction $W=((u_1, \dots, u_k,v),(\mu(1), \dots, \mu(k)))$ of $W_0$ can be effectively calculated from the above and the induction hypothesis, and
\[
m^*_{\T,P}(w) = u_1 \cdots u_kv\,.
\]
\qed

\bigskip\noindent
{\bf Corollary 4.16.}
{\it \begin{itemize}
\item[(1)]
$G(M)$ has a solution to the word problem.
\item[(2)]
There exists an algorithm which, given $w \in \AA(M)^{\pm *}$, decides whether $\overline{w} \in M$.
\item[(3)]
Let $P$ be a parabolic submonoid of $M$. There exists an algorithm which, given $w \in \AA(M)^{\pm *}$, decides whether $\overline{w} \in G(P)$.
\end{itemize}}

\bigskip\noindent
{\bf Proof.}
As pointed out in Corollary 4.13, a FC tree $\T$ for $M$ can be effectively calculated. Let $w \in \AA(M)^{\pm *}$. Then we have $\overline{w} = 1$ if and only if $m_{\T,\{1\}}^*(w)=1$, and we have $\overline{w} \in M$ if and only if $m_{\T,\{1\}}^*(w)\in \SS(M)^*$. Let $P$ be a parabolic submonoid of $M$. Then $\overline{w} \in G(P)$ if and only if $m_{\T,P}^*(w)=1$.
\qed



\bigskip\bigskip\noindent
{\bf Eddy  Godelle,}

\smallskip\noindent
Université de Caen, Laboratoire LMNO, UMR 6139 du CNRS, Campus II, 14032 Caen cedex, France. 

\smallskip\noindent
E-mail: {\tt eddy.godelle@unicaen.fr}

\bigskip\noindent
{\bf Luis Paris,}

\smallskip\noindent 
Université de Bourgogne, Institut de Mathématiques de Bourgogne, UMR 5584 du CNRS, B.P. 47870, 21078 Dijon cedex, France.

\smallskip\noindent
E-mail: {\tt lparis@u-bourgogne.fr}



\begin{thebibliography}{99}

\small

\bibitem{Alt}
{\bf J.\,A. Altobelli.}
{\it The word problem for Artin groups of FC type.}
J. Pure Appl. Algebra {\bf 129} (1998), no. 1, 1--22.

\bibitem{AltCha1}
{\bf J. Altobelli, R. Charney.}
{\it A geometric rational form for Artin groups of FC type.}
Geom. Dedicata {\bf 79} (2000), no. 3, 277--289.

\bibitem{Arnol1}
{\bf  V.\,I. Arnol'd.}
{\it Braids of algebraic functions and cohomologies of swallowtails. (Russian)}
Uspehi Mat. Nauk {\bf 23} (1968), no. 4, 247--248. 

\bibitem{Arnol2}
{\bf V.\,I. Arnol'd.}
{\it The cohomology classes of algebraic functions that are preserved under Tschirnhausen transformations. (Russian)}
Funkcional. Anal. i Prilo\v zen. {\bf 4} (1970), no. 1, 84--85.

\bibitem{Arnol3}
{\bf V.\,I. Arnol'd.}
{\it Certain topological invariants of algebrac functions. (Russian)}
Trudy Moskov. Mat. Ob\v s\v c. {\bf 21} (1970), 27--46. 

\bibitem{Arnol4}
{\bf V.\,I. Arnol'd.}
{\it Topological invariants of algebraic functions. II. (Russian)}
Funkcional. Anal. i Prilo\v zen. {\bf 4} (1970), no. 2, 1--9.
Translation in Functional Anal. Appl. {\bf 4} (1970), 91--98.

\bibitem{BiMaMe}
{\bf J.-C. Birget, S.\,W. Margolis, J. Meakin.}
{\it On the word problem for tensor products and amalgams of monoids.} 
Internat. J. Algebra Comput. {\bf 9} (1999), no. 3-4, 271--294.

\bibitem{Bries1}
{\bf E. Brieskorn.}
{\it Singular elements of semi-simple algebraic groups.}
Actes du Congrès International des Mathématiciens (Nice, 1970), Tome 2, pp. 279--284. 
Gauthier-Villars, Paris, 1971. 

\bibitem{Bries2}
{\bf E. Brieskorn.}
{\it Die Fundamentalgruppe des Raumes der regulären Orbits einer endlichen komplexen Spiegelungsgruppe.}
Invent. Math. {\bf 12} (1971), 57--61. 

\bibitem{Bries3}
{\bf E. Brieskorn.}
{\it Sur les groupes de tresses [d'après V. I. Arnol'd].}
Séminaire Bourbaki, 24ème année (1971/1972), Exp. No. 401, pp. 21--44. 
Lecture Notes in Math., Vol. 317, Springer, Berlin, 1973. 

\bibitem{BriSai1}
{\bf E. Brieskorn, K. Saito.}
{\it Artin-Gruppen und Coxeter-Gruppen.}
Invent. Math. {\bf 17} (1972), 245--271.

\bibitem{Charn1}
{\bf R. Charney.}
{\it Artin groups of finite type are biautomatic.}
Math. Ann. {\bf 292} (1992), no. 4, 671--683. 

\bibitem{Charn2}
{\bf R. Charney.}
{\it Geodesic automation and growth functions for Artin groups of finite type.}
Math. Ann. {\bf 301} (1995), no. 2, 307--324.

\bibitem{ChaDav1}
{\bf R. Charney, M.\,W. Davis.,}
{\it The $K(\pi,1)$-problem for hyperplane complements associated to infinite reflection groups.}
J. Amer. Math. Soc. {\bf 8} (1995), no. 3, 597--627. 

\bibitem{ChaDav2}
{\bf R. Charney, M.\,W. Davis.}
{\it Finite $K(\pi,1)$'s for Artin groups.}
Prospects in topology (Princeton, NJ, 1994), 110--124, Ann. of Math. Stud., 138, Princeton Univ. Press, Princeton, NJ, 1995. 

\bibitem{Cherma1}
{\bf A. Chermak.}
{\it Locally non-spherical Artin groups.}
J. Algebra {\bf 200} (1998), no. 1, 56--98.

\bibitem{Deh2}
{\bf P. Dehornoy.}
{\it Gaussian groups are torsion free.}
J. Algebra {\bf 210} (1998), 291--297.

\bibitem{Deh}
{\bf P. Dehornoy.}
{\it Groupes de Garside.}
Ann. Sci. \'Ecole Norm. Sup. (4) {\bf 35} (2002), no. 2, 267--306.

\bibitem{Deh3}
{\bf P. Dehornoy.}
{\it Alternating normal forms for braids and locally Garside monoids.}
J. Pure Appl. Algebra {\bf 212} (2008), no. 11, 2413--2439. 

\bibitem{Deh4}
{\bf P. Dehornoy.}
{\it Left-Garside categories, self-distributivity, and braids.}
Ann. Math. Blaise Pascal {\bf 16} (2009), no. 2, 189--244. 

\bibitem{DDGKM1}
{\bf P. Dehornoy, F. Digne, E. Godelle, D. Krammer, J. Michel.}
{\it Garside theory.}
Book in preparation. Avaible at http://www.math.unicaen.fr/~garside/Garside.pdf

\bibitem{DePa}
{\bf P. Dehornoy, L. Paris.}
{\it Gaussian groups and Garside groups, two generalisations of Artin groups.}
Proc. London Math Soc. (3) {\bf 79} (1999), 569--604.

\bibitem{Delig1}
{\bf P. Deligne.}
{\it Les immeubles des groupes de tresses g\'en\'eralis\'ees.} 
Invent. Math. {\bf 17} (1972), 273--302. 

\bibitem{Digne1} 
{\bf F. Digne.}
{\it Présentations duales des groupes de tresses de type affine $\widetilde A$.}
Comment. Math. Helv. {\bf 81} (2006), no. 1, 23--47.

\bibitem{Digne2}
{\bf F. Digne.}
{\it A Garside presentation for Artin-Tits groups of type $\tilde C_n$.}
Preprint, arXiv: 1002.4320. 

\bibitem{DiM}
{\bf F. Digne, J. Michel.}
{\it Garside and locally Garside categories.}
Preprint, arXiv: math/0612652. 

\bibitem{Eps}
{\bf D.\,B.\,A. Epstein.}
{\it A result on free products with amalgamation.}
J. London Math. Soc. {\bf 37} (1962), 130--132.

\bibitem{ECHLPT1}
{\bf D.\,B.\,A. Epstein, J.\,W. Cannon, D.\,F. Holt, S.\,V.\,F. Levy, M.\,S. Paterson, W.\,P. Thurston.}
{\it Word processing in groups.}
Jones and Bartlett Publishers, Boston, MA, 1992.

\bibitem{Garsi1}
{\bf F.\,A. Garside.}
{\it The braid group and other groups.}
Quart. J. Math. Oxford Ser. (2) {\bf 20} (1969), 235--254. 

\bibitem{GoG}
{\bf V. Gebhardt, J. Gonz\'alez-Meneses.}
{\it Solving the conjugacy problem in Garside groups by cyclic sliding.}
J. Symbolic Comput. {\bf 45} (2010), no. 6, 629--656.

\bibitem{Godel1}
{\bf E. Godelle.}
{\it Parabolic subgroups of Artin groups of type FC.}
Pacific J. Math. {\bf 208} (2003), no. 2, 243--254.

\bibitem{Godel2}
{\bf E. Godelle.}
{\it Artin-Tits groups with CAT(0) Deligne complex.}
J. Pure Appl. Algebra {\bf 208} (2007), no. 1, 39--52. 

\bibitem{God}
{\bf E. Godelle.}
{\it Parabolic subgroups of Garside groups.}
J. of Algebra {\bf 317} (2007), 1--16.

\bibitem{Godel3}
{\bf E. Godelle.}
{\it Parabolic subgroups of Garside groups II: ribbons.}
J. Pure Appl. Algebra {\bf 214} (2010), no. 11, 2044--2062. 


\bibitem{GodPar1}
{\bf E. Godelle, L. Paris.}
{\it $K(\pi,1)$ and word problems for infinite type Artin-Tits groups, and applications to virtual braid groups.}
Math. Z., to appear.

\bibitem{Gorju1}
{\bf V.\,V. Gorjunov.}
{\it The cohomology of braid groups of series C and D and certain stratifications.}
Funktsional. Anal. i Prilozhen. {\bf 12} (1978), no. 2, 76--77. 
Translation in: Functional Anal. Appl. {\bf 12} (1978), no. 2, 139--149.

\bibitem{Gre}
{\bf J. Green.}
{\it On the structure of semigroups.} 
Ann. of Math. (2) {\bf 54} (1951), 163--172.

\bibitem{How2}
{\bf J.\,M. Howie.}
{\it Embedding theorems with amalgamation for semigroups.} 
Proc. London Math. Soc. (3) {\bf 12} (1962), 511--534.

\bibitem{How}
{\bf J.\,M. Howie.}
{\it An introduction to semigroup theory.} 
L.M.S. Monographs, No. 7. Academic Press, London -- New York, 1976.

\bibitem{Lin1}
{\bf V.\,Y. Lin.}
{\it Artin braids and the groups and spaces connected with them.}
Itogi Nauki i Tekhniki, Algebra, Topologiya, Geometriya 17, VINITI, Moscow, 1979, pp. 159--227. 
Translation in: Journal of Soviet Math. {\bf 18} (1982), 736--788.

\bibitem{Mic}
{\bf J. Michel.}
{\it A note on words in braid monoids.}
J. Algebra {\bf 215} (1999), 366--377.

\bibitem{Paris1}
{\bf L. Paris.}
{\it Universal cover of Salvetti's complex and topology of simplicial arrangements of hyperplanes.}
Trans. Amer. Math. Soc. {\bf 340} (1993), no. 1, 149--178.

\bibitem{Par}
{\bf L. Paris.}
{\it Artin monoids inject in their groups.}
Comment. Math. Helv. {\bf 77} (2002), 609--637.

\bibitem{Sch}
{\bf O. Schreier.}
{\it Die Untergruppen der Freien Gruppen.}
Abh. Math. Sem. Univ. of Hamb. {\bf 5} (1927), 161--183.

\bibitem{Ser}
{\bf J.-P. Serre.}
{\it Arbres, amalgames, $SL_2$.}
Ast\'erisque {\bf 46}, SMF, Paris, 1977.

\bibitem{VdL}
{\bf H. Van der Lek.}
{\it The homotopy type of complex hyperplane complements.}
Ph. D. thesis, Nijmegen, 1983.

\end{thebibliography}
\end{document}